\documentclass[11pt]{amsart}
\usepackage[a4paper,margin=27mm]{geometry}
\usepackage{amssymb,mathtools}
\usepackage[T1]{fontenc}
\usepackage{lmodern}
\usepackage[colorlinks=true,linkcolor=blue,citecolor=blue,urlcolor=blue]{hyperref}
\newtheorem{theorem}{Theorem}[section]
\newtheorem{proposition}[theorem]{Proposition}
\newtheorem{lemma}[theorem]{Lemma}
\theoremstyle{remark}
\newtheorem{remark}[theorem]{Remark}
\numberwithin{equation}{section}
\DeclareMathOperator{\Res}{Res}
\DeclareMathOperator{\Tr}{Tr}
\DeclareMathOperator{\ord}{ord}
\newcommand{\eps}{\varepsilon}
\newcommand{\Q}{\mathbb Q}
\newcommand{\Z}{\mathbb Z}

\title[A secondary term in the cubic first moment]{A secondary term in the first moment of
primitive cubic Dirichlet \(L\)-functions}
\author{Ziwei Hong}
\address{School of Mathematics, Renmin University of China}
\address{School of Mathematics and Statistics, University of New South Wales}
\email{hongziwei@live.com}
\date{24 October 2026}

\begin{document}

\begin{abstract}
P. Gao and L. Zhao obtained an asymptotic formula for the smoothed first
moment of primitive cubic Dirichlet \(L\)-functions over \(\Q\) with
conductors coprime to \(3\), with a main
term of order \(Q\) and error \(O\left(Q^{\frac{7}{8}+\eps}\right)\)
under the generalized Lindel\"of hypothesis. Building on their double
Dirichlet series, we assume the generalized Riemann hypothesis and
separate an explicit secondary term of order \(Q^{\frac{5}{6}}\)
in the same family, with error \(O\left(Q^{\frac{23}{28}+\eps}\right)\).
We explicitly evaluate the secondary coefficient and prove that it
is strictly negative. We also prove unconditionally that the secondary
polynomial in the squarefree second moment of the natural cubic Hecke
\(L\)-functions over \(\Q(\sqrt{-3})\) has degree one and strictly
positive leading coefficient, which we evaluate explicitly.
\end{abstract}
\maketitle

\section{Introduction}

We consider the smoothed first moment of primitive cubic Dirichlet
\(L\)-functions over \(\Q\). More precisely, we study
\[
\sum_{(q,3)=1}\sideset{}{^*}\sum_{\substack{\chi\bmod q\\
\chi^3=\chi_0}}L\left(\frac{1}{2},\chi\right)\Phi\left(\frac{q}{Q}\right),
\]
where \(Q>0\) tends to infinity and \(q\) runs over positive integers
coprime to \(3\). The inner sum is over Dirichlet characters modulo
\(q\), and \(\chi_0\) denotes the principal character modulo \(q\).
The star restricts the sum to primitive characters of exact order
\(3\), so \(q\) is their conductor and the character of conductor
\(1\) is excluded. Here \(L(s,\chi)\) denotes the Dirichlet
\(L\)-function attached to \(\chi\), evaluated at its central point
in the sum above. The weight \(\Phi(x)\) is a fixed nonnegative smooth
function compactly supported on the set of positive real numbers
\(\mathbb R_+\), so \(\Phi\left(\frac{q}{Q}\right)\) restricts the
conductors to size \(Q\). We write \(\widehat\Phi\) for its Mellin transform,
\[
\widehat\Phi(s)=\int_0^\infty\Phi(x)x^{s-1}\,dx.
\]

P. Gao and L. Zhao \cite{GZ} study the first moment through a double Dirichlet
series \(A_3(s,w)\). Under the generalized Lindel\"of hypothesis, their
method gives a main term of order \(Q\) with error
\(O\left(Q^{\frac{7}{8}+\eps}\right)\). Their continuation also exhibits the polar
hyperplane \(s+w=\frac{4}{3}\), whose intersection with \(w=\frac{1}{2}\) corresponds
to a term of order \(Q^{\frac{5}{6}}\). To separate this contribution, one
needs both continuation to the left of \(s=\frac{5}{6}\) and suitable bounds
there.

More recently, P. Gao and L. Zhao \cite{GZquadratic} applied the method
of M. \v Cech \cite{Cechfirst} to the first moment of primitive quadratic
Hecke \(L\)-functions over \(\Q(i)\). Under GRH, they obtained a secondary
term of the form \(X^{\frac{1}{3}}\) times a linear polynomial in
\(\log X\), with error \(O\left(X^{\frac{1}{4}+\eps}\right)\).
Their argument uses the explicit evaluation of quadratic Gauss sums
and the fact that the square of a quadratic character is principal.
Consequently, the \(L\)-functions attached to these squared characters
reduce to the Dedekind zeta function with finitely many Euler factors
removed. This simplification does not hold for a nontrivial cubic
character, whose square is still nonprincipal. Thus this part of
their method does not extend directly to the cubic family considered here.

A. Diaconu, B. Ion, V. Pa\c{s}ol, and A. A. Popa \cite{DIPP} also
separate a secondary term in the first moment of cubic Hecke
\(L\)-functions. Their results assume that the base field contains
the required roots of unity; in the cubic case, it contains a primitive
cube root of unity. This is the Kummer setting.

M. H. Hamdar \cite{Hamdar} studies the first moment of cubic Hecke
\(L\)-functions over \(\Q(\sqrt{-3})\), with squarefree parameters
congruent to \(1\) modulo \(9\). At the central point, he obtains
unconditionally a main term of order \(Q\) and a secondary term of
order \(Q^{\frac{5}{6}}\), with error
\(O\left(Q^{\frac{4}{5}+\eps}\right)\).
He also studies values away from the central point and identifies a
transition at \(s=\frac{1}{3}\). His proof uses the approximate
functional equation, the cubic large sieve, and explicit residues of
Dirichlet series associated with cubic Gauss sums. A cancellation
between contributions from the principal and dual sums permits the
secondary term to be separated.
The secondary coefficient is expressed as a residue, but its
nonvanishing at the central point is not explicitly verified
in \cite[Section 4.2]{Hamdar}.

A. de Faveri \cite[Theorem 1.3 and footnote 5]{deFaveri} obtains an
unconditional second moment asymptotic for cubic Hecke \(L\)-functions
over squarefree ideals in number fields containing a primitive cube
root of unity. His formula includes a secondary contribution
\(X^{\frac{5}{6}}Q_W(\log X)\), with error
\(O\left(X^{\frac{3}{4}+\eps}\right)\), but does not rule out
\(Q_W\equiv0\). Theorem~\ref{thm:second-moment} below proves that,
for the natural family over \(\Q(\sqrt{-3})\), this polynomial has
degree exactly one and positive leading coefficient.
For the first moment of primitive cubic Dirichlet
\(L\)-functions over \(\Q\) considered here, we explicitly evaluate
the secondary coefficient and prove that it is strictly negative.
Thus our asymptotic formula under GRH has a nonzero secondary term.

The results of \cite{DIPP,Hamdar} concern Hecke \(L\)-functions over
fields containing a primitive cube root of unity. Our first moment result
concerns primitive cubic Dirichlet \(L\)-functions over \(\Q\), which
does not contain a primitive cube root of unity. This result is therefore
in the non-Kummer setting. Although cubic residue symbols over
\(\Q(\sqrt{-3})\) also parameterize our characters, the associated
Dirichlet \(L\)-functions differ from the Hecke \(L\)-functions
studied in those works. Moreover, the parameterization of primitive
Dirichlet characters requires the squarefree parameters to have no
rational prime divisor. The theorems of \cite{DIPP,Hamdar} therefore
do not directly give an asymptotic formula for our family.
Our first contribution is to separate the secondary term in the family
studied by Gao and Zhao \cite{GZ}, under GRH, and to evaluate its
coefficient explicitly and prove that it is nonzero.
Diaconu, Ion, Pa\c{s}ol, and Popa \cite[Remark 1.4.1]{DIPP} observe
that their assumption on roots of unity can in principle be relaxed
with additional technical work; in this direction, the joint estimate
in Lemma~\ref{lem:blocks} provides the analytic input needed to separate
the secondary term for the cubic first moment over \(\Q\) under GRH.

Throughout the paper,
\[
K=\Q(\sqrt{-3}),\qquad
\omega=\frac{-1+\sqrt{-3}}2,\qquad
c_K=\sqrt{D_K}=\sqrt{-3},\qquad
r_K=\Res_{s=1}\zeta_K(s)=\frac{\pi}{3\sqrt3}.
\]
A primary element of \(\Z[\omega]\) is one congruent to \(1\) modulo
\(3\). For a rational prime \(p\ne3\), put
\begin{equation}
\kappa_p=\prod_{\varpi\mid p}\left(1+\frac1{N(\varpi)}\right)^{-1}
=\begin{cases}
(1+p^{-1})^{-2},&p\equiv1\pmod3,\\
(1+p^{-2})^{-1},&p\equiv2\pmod3,
\end{cases}
\label{eq:kappap}
\end{equation}
where \(\varpi\) runs over the primary prime factors in \(K\).
This factor occurs in both first-moment residues.

\begin{theorem}\label{thm:main}
Assume GRH for Dirichlet \(L\)-functions and Hecke \(L\)-functions
associated with characters of finite order over \(K\). With \(\Phi\)
as above, for every \(\eps>0\), as \(Q\to\infty\),
\begin{align}
\sum_{(q,3)=1}\sideset{}{^*}\sum_{\substack{\chi\bmod q\\
\chi^3=\chi_0}}L\left(\frac{1}{2},\chi\right)\Phi\left(\frac{q}{Q}\right)
={}&C_3\widehat\Phi(1)Q
+C_3'\widehat\Phi\left(\frac{5}{6}\right)Q^{\frac{5}{6}}
+O_{\Phi,\eps}\left(Q^{\frac{23}{28}+\eps}\right),                 \label{eq:main}
\end{align}
where \(C_3\) is the constant in the main term of \cite{GZ} at the
central point. Explicitly,
\begin{equation}
C_3=\frac{3r_K\zeta\left(\frac{3}{2}\right)}{4\zeta_K(2)}
\prod_{p\ne3}
\left(1-\kappa_pp^{-2}
+\bigl(\kappa_p(1+p^{-2})-1\bigr)p^{-\frac{3}{2}}\right).
\label{eq:C3}
\end{equation}
The secondary coefficient is
\begin{align}
C_3'={}&
\frac{(2\pi)^{\frac{5}{3}}\zeta\left(\frac{3}{2}\right)\zeta\left(\frac{5}{6}\right)}
 {8\,3^{\frac{3}{2}}\Gamma\left(\frac{2}{3}\right)\zeta_K(2)}
\left(1-3^{-\frac{5}{3}}\right)                                      \notag\\
&\quad\times
\prod_{p\equiv1\ (3)}
\frac{\left(1-p^{-\frac{5}{6}}\right)\left(1+2p^{-1}+p^{-\frac{5}{6}}-2p^{-\frac{5}{2}}\right)}
 {(1+p^{-1})^2}
\prod_{p\equiv2\ (3)}
\frac{1-p^{-\frac{5}{3}}}{1+p^{-2}}.
\label{eq:C3prime}
\end{align}
Both Euler products in \eqref{eq:C3prime} converge absolutely, and
\(C_3'<0\).
\end{theorem}

Since \(\frac{5}{6}-\frac{23}{28}=\frac{1}{84}\), choosing \(0<\eps<\frac{1}{84}\) separates the
secondary contribution. If \(\Phi\) is not identically zero, then
\(\widehat\Phi\left(\frac{5}{6}\right)>0\), so this contribution is nonzero.

In the expression for \(C_{3,1}(s,w)\) in \cite{GZ}, M\"obius inversion introduces
a sum over \(d\) with coefficient \(\mu_{\Q}(d)\). We estimate this sum
together with the sum over \(m\), retaining cancellation in both
variables. The joint estimate in Lemma~\ref{lem:blocks} is obtained
from a ratio of Dirichlet \(L\)-functions under GRH. We combine it
with the bound for Dirichlet series associated with Gauss sums in
C. David, A. de Faveri, A. Dunn, and J. Stucky \cite[(6.13)]{DDFDS}.
Interpolation and Bochner's tube theorem then give analytic
continuation to \(\Re s>\frac{23}{28}\) when \(w=\frac{1}{2}\),
which allows the Mellin contour to pass the pole at \(s=\frac{5}{6}\).
GRH is used both in the estimate for the ratio and in the estimates
for \(L\)-values that give the initial domains of continuation.
The formulas for the Gauss sum residues are unconditional.

We retain the notation \(A_3,C_{3,1},H_{3,1}\), the residue symbols,
and the Gauss sums of \cite{GZ} for the first moment.
The Hecke family in \cite{DDFDS} supplies analytic inputs for that
argument. Our second result below concerns the second moment of
the natural Hecke family itself.

Our second result is unconditional and concerns the natural family of cubic Hecke
\(L\)-functions over \(K\). Put \(\lambda=c_K\), let \(S\) consist of
the infinite place and the prime \((\lambda)\), and write \(I(S)\) for
the nonzero integral ideals coprime to \(3\). Each ideal in \(I(S)\)
has a unique primary generator. For such a generator \(a\), let
\(\chi_a\) be the Hecke character determined by
\[
\chi_a((b))=\left(\frac{a}{b}\right)_3
\quad\text{for primary }b\text{ with }(b,3a)=1.
\]
We write \(L^S\) for the corresponding Hecke \(L\)-function with its
Euler factor at \(\lambda\) removed, and \(\mu_K\) for the ideal
M\"obius function. For a fixed nonnegative, nonzero
\(W\in C_c^\infty((0,\infty))\), define
\begin{equation}
\mathcal M_W(X)=
\sum_{\mathfrak a\in I(S)}\mu_K^2(\mathfrak a)
\left|L^S\left(\frac12,\chi_{\mathfrak a}\right)\right|^2
W\left(\frac{N\mathfrak a}{X}\right),
\qquad
\widehat W(s)=\int_0^\infty W(x)x^{s-1}\,dx.
\label{eq:second-moment-family}
\end{equation}
Here \(\chi_{\mathfrak a}\) denotes \(\chi_a\) for the primary generator
of \(\mathfrak a\). The parameter set contains all squarefree ideals
in \(I(S)\); in particular, it has neither the rational prime divisor
restriction of Theorem~\ref{thm:main} nor a fixed congruence class
modulo \(9\).

\begin{theorem}\label{thm:second-moment}
With \(W\) as above, there are a real polynomial \(P_W\) of degree at most
one and a real constant \(b_W\) such that, for every \(\eps>0\),
\begin{equation}
\begin{split}
\mathcal M_W(X)
={}&XP_W(\log X)\\
&+X^{\frac56}
\left(A_{\mathrm{sf}}\widehat W\left(\frac56\right)\log X+b_W\right)
+O_{W,\eps}\left(X^{\frac34+\eps}\right)
\end{split}
\label{eq:second-moment-asymptotic}
\end{equation}
as \(X\to\infty\), where
\begin{equation}
\begin{split}
A_{\mathrm{sf}}
={}&
\frac{2(2\pi)^4(1+3^{-\frac12})\Gamma\left(\frac16\right)}
{3^8\Gamma\left(\frac56\right)\Gamma\left(\frac23\right)^2}
\zeta_\lambda\left(\frac12\right)^2\\
&\qquad\times
\prod_{p\notin S}
\left[
(1-q^{-1})^2
\left(1+2q^{-1}-3q^{-2}+2q^{-\frac52}\right)
\right]>0.
\end{split}
\label{eq:second-coefficient}
\end{equation}
Here \(p\) runs over prime ideals of \(K\), \(q=Np\), and
\(\zeta_\lambda(s)=(1-3^{-s})\zeta_K(s)\).
The product in \eqref{eq:second-coefficient} converges absolutely
when each bracket is regarded as one local factor.
Consequently, the secondary polynomial in
\eqref{eq:second-moment-asymptotic} has degree exactly one and
strictly positive leading coefficient.
\end{theorem}

The asymptotic expansion and its error term are supplied by
\cite[Theorem~1.3]{deFaveri}. Our contribution to this second moment
is the explicit formula \eqref{eq:second-coefficient} and its strict
positivity for the natural family over \(\Q(\omega)\).
This settles the possible vanishing of the secondary polynomial
in that specialization. The constant \(b_W\) is not needed for this
conclusion, and we make no assertion about its sign. Thus the
positive coefficient in Theorem~\ref{thm:second-moment} and the
negative coefficient \(C'_3\) in Theorem~\ref{thm:main} concern
different moments and different families.

The proof is given in Section~\ref{sec:second-moment}. Patterson's
residue formula, in the form recorded in \cite{DDFDS}, expresses
the square sum of the raw cubic residues as a quotient of zeta
functions. This quotient has a simple zero at \(u=\frac16\), with
positive derivative. We compute the complete finite pairing in
the multiple Dirichlet series of \cite{DIPP}, including the
primitive root numbers and the actual Euler factor at \(\lambda\),
and then evaluate the local factors of the squarefree sieve.
The resulting factors are positive. Finally, the estimates of
\cite{DIPP,deFaveri} give normal convergence after the poles are
removed, so the sieve can be interchanged with extraction of the
Laurent coefficient. This argument uses no GRH input from the
proof of Theorem~\ref{thm:main}.

Section~\ref{sec:preliminaries} records the analytic and arithmetic inputs,
and Section~\ref{sec:first-moment-proof} establishes the continuation and
contour shift for the first moment. The evaluations of both secondary
terms are collected in Section~\ref{sec:secondary-terms}.
Section~\ref{sec:notation-comparison} first explains the relation between
our notation and that of \cite{deFaveri}; Sections~\ref{sec:first-secondary-residue}
and~\ref{sec:second-moment} then give the first- and second-moment
calculations, respectively.

\section{Preliminaries}\label{sec:preliminaries}

\subsection{Residue symbols}

For a primary \(n\), let
\(\chi_{3,n}(x)=\left(\frac{x}{n}\right)_3\) and
\(\widehat\chi_{3,n}(m)=\left(\frac{m}{n}\right)_3\) for rational integers \(m\).
When \(n\) is squarefree and has no rational prime divisor,
\(\widehat\chi_{3,n}\) is primitive of conductor \(N(n)\).
For arbitrary squarefree primary \(n\), the same notation on rational
integers denotes a possibly imprimitive Dirichlet character.

Let \(\mu_{\Q}(d)=\mu(|d|)\), where \(\mu\) is the usual M\"obius
function.

\subsection{Gauss sums}

We use precisely the additive normalization of \cite{GZ}:
\begin{equation}
g_K(r,\chi_{3,n})=
\sum_{x\bmod n}\chi_{3,n}(x)
 \exp\left(2\pi i\Tr\left(\frac{rx}{c_Kn}\right)\right),
\qquad g_K(\chi_{3,n})=g_K(1,\chi_{3,n}).
\label{eq:gauss}
\end{equation}
Here \(\Tr\) denotes the trace from \(K\) to \(\Q\).
Thus the cubic Gauss sum normalized using the trace in \cite{DDFDS} is
\(g_K(c_Kr,\chi_{3,n})\). This change of additive normalization is
necessary when importing their residue formulas.

If \(d\ge1\) is squarefree and \((d,3)=1\), take its primary associate.
Then
\begin{equation}
g_K(c_K,\chi_{3,d})=d,\qquad
\left(\frac{m}{d}\right)_3=1\quad\text{when }(m,d)=1.
\label{eq:rational-gauss}
\end{equation}
For an inert prime, summing the cubic character over trace fibers in
\(\mathbb F_{p^2}\) gives \(p\). For a split prime, conjugation pairs
the two Gauss sums normalized using the trace, whose product is \(p\);
the two reciprocity factors multiply to \(1\). Cubic reciprocity extends the identity
to squarefree rational \(d\). Changing the sign to its primary
associate has no effect. Here and below powers of a rational primary
associate mean powers of its positive absolute value.

\subsection{Estimating \texorpdfstring{\(L\)-functions}{L-functions}}

We use the standard consequences of GRH that are uniform in the conductor. For fixed
\(\delta,\eta>0\), to the right of \(\frac{1}{2}+\delta\),
\[
L(z,\chi),\quad L(z,\chi)^{-1}
\ll_{\delta,\eta}\bigl(q(1+|\Im z|)\bigr)^\eta,
\]
for a primitive nonprincipal Dirichlet character of conductor \(q\),
with real parts restricted to a fixed bounded interval. For a
principal numerator, its pole at \(z=1\) must be removed first.
The reciprocal estimate is the analogue for general characters of
\cite[(2.21)]{Cech}, which is stated there for quadratic characters.
We use the analogous Lindel\"of consequence for the Hecke
\(L\)-functions occurring in the region obtained from the original series.

For imprimitive characters, the additional Euler factors and their
reciprocals cost at most an arbitrarily small power of the modulus.
Indeed, for fixed \(\delta>0\),
\[
\prod_{p\mid q}\left(1-p^{-\frac{1}{2}-\delta}\right)^{-1}
\ll_{\delta,\eta}q^\eta.
\]
The parameters \(\delta\) and \(\eta\) can be chosen independently.

\subsection{Analytic behavior of Dirichlet series associated with Gauss sums}

Following \cite{GZ}, define the Dirichlet series associated with Gauss sums
\begin{equation}
H_{3,1}(s,l)=
\sum_{\substack{n\ {\rm primary}\\(n,l)=1}}
\frac{g_K(l,\chi_{3,n})}{N(n)^s}.
\label{eq:H31}
\end{equation}
For rational \(r\ge1\), the series \(H_{3,1}(s,rc_K)\) is the
coprime series normalized using the trace in \cite{DDFDS}.

Put
\[
\zeta_{c_K}(s)=(1-3^{-s})\zeta_K(s).
\]
This is the function denoted \(\zeta_\lambda\) in \cite{DDFDS}, where
\(\lambda=c_K\).

\begin{lemma}\label{lem:Hbound}
For every rational integer \(r\ge1\),
\[
\left(s-\frac{4}{3}\right)\zeta_{c_K}(3s-2)H_{3,1}(s,rc_K)
\]
is holomorphic in \(\Re s>\frac{2}{3}\). For every compact
\(I\subset\left(\frac{2}{3},\frac{3}{2}\right)\) and \(\eps>0\), it satisfies
\begin{equation}
\left|\left(s-\frac{4}{3}\right)\zeta_{c_K}(3s-2)H_{3,1}(s,rc_K)\right|
\ll_{I,\eps}r^{\frac{3}{2}-\Re s+\eps}(1+|\Im s|)^A,
\qquad \Re s\in I,
\label{eq:Hbound}
\end{equation}
where \(A\) depends only on \(I\) and \(\eps\).
\end{lemma}

\begin{proof}
The unrestricted series \(\psi(r,s)\) in \cite{DDFDS} has coefficients
\(g_K(c_Kr,\chi_{3,n})\). Their equation (6.13), including the
holomorphy statement immediately preceding it, gives
\begin{equation}
\psi(r,s)\zeta_{c_K}(3s-2)
\ll_{\sigma_0}
\frac{|s|}{|s-\frac{4}{3}|}
 N(r)^{\frac{\sigma_0-\Re s}{2}}
 (1+|\Im s|)^{2(\sigma_0-\Re s)}
\label{eq:completed-input}
\end{equation}
for \(\frac{3}{2}<\sigma_0<2\) and
\(2-\sigma_0\le\Re s\le\sigma_0\), with the pole at \(\frac{4}{3}\)
removed when discussing holomorphy. We keep the zeta factor in this
estimate.

To insert coprimality, take primary divisors and factor the twist in
\(\Z[\omega]\) as
\[
r=\eta ab^2c^3,
\]
where \(a,b\) are squarefree and coprime, and \(\eta\) belongs to a
fixed finite set of ramified unit twists. Cubes can be removed from
the twist in the coprime series. Apply \cite[Corollary 6.2]{DDFDS}
to remove the conditions at the primes dividing \(abc\).
Using the square root bounds for the Gauss sums in that formula and
taking \(\sigma_0\) arbitrarily close to \(\frac{3}{2}\), the resulting bound,
apart from \(N(r)^\eps(1+|\Im s|)^A\), is
\begin{equation}
N(ab^2)^{\frac{3}{4}-\frac{\sigma}{2}}
\sum_{d\mid a}N(d)^{\frac{3}{4}-\frac{3\sigma}{2}}
\sum_{\substack{e\mid c\\ e\ {\rm squarefree}\\(e,ab)=1}}
N(e)^{\frac{5}{4}-\frac{3\sigma}{2}},
\qquad \sigma=\Re s.
\label{eq:divisor-bound}
\end{equation}
The local denominators are products of
\(1-N(\varpi)^{2-3s}\). They have no zeros for \(\sigma>\frac{2}{3}\),
and their reciprocals are \(\ll N(abc)^\eps\) on every fixed smaller
strip.

The first divisor exponent in \eqref{eq:divisor-bound} is negative.
If the second is positive, then
\[
\frac{5}{4}-\frac{3\sigma}{2}
\le 3\left(\frac{3}{4}-\frac{\sigma}{2}\right).
\]
Thus the second sum is absorbed by the cube part, with a divisor
loss \(N(c)^\eps\). The result is
\(N(r)^{\frac{3}{4}-\frac{\sigma}{2}+\eps}\). Since \(N(r)=r^2\) for rational
\(r\), this gives \eqref{eq:Hbound}. The same finite coprimality
formula proves holomorphy; its denominators have zeros only on
\(\Re s=\frac{2}{3}\). Absolute convergence covers \(\Re s>\frac{3}{2}\).
\end{proof}

We also need the residue of \(H_{3,1}(s,rc_K)\).

\begin{lemma}\label{lem:Hresidue}
Write a positive rational integer uniquely as
\[
r=3^\nu ab^2c^3,\qquad
\nu\in\{0,1,2\},\quad
\mu_{\Q}^2(ab)=1,\quad (a,b)=1,\quad (ab,3)=1.
\]
The integer \(c\) is unrestricted and may share primes with \(3ab\).
Then
\begin{equation}
\Res_{s=\frac{4}{3}}H_{3,1}(s,rc_K)
=\begin{cases}
\displaystyle
\frac{(2\pi)^{\frac{5}{3}}}
 {8\,3^{\frac{3}{2}}\Gamma\left(\frac{2}{3}\right)\zeta_K(2)}
 (3^\nu a)^{-\frac{1}{3}}
 \prod_{\substack{p\mid ac\\p\ne3}}\kappa_p,
 &b=1,\ \nu=0,1,\\[8pt]
0,&\text{otherwise}.
\end{cases}
\label{eq:Hresidue}
\end{equation}
\end{lemma}

\begin{proof}
Use the correspondence between normalizations after \eqref{eq:gauss}.
The unrestricted residue in \cite[Lemma 6.3]{DDFDS} is
\(\frac{c_0\tau_3(r)}{N(r)^{\frac{1}{6}}}\), where their equation (6.10) gives
\[
c_0=\frac{(2\pi)^{\frac{5}{3}}}
 {8\,3^{\frac{9}{2}}\Gamma\left(\frac{2}{3}\right)\zeta_K(2)}.
\]
Their theta coefficient formula (6.8) gives
\(\tau_3(1)=\tau_3(3)=27\) and \(\tau_3(9)=0\).
After division by \(N(r)^{\frac{1}{6}}\), the relative ramified factors
are \(1,3^{-\frac{1}{3}},0\). Cubic invariance removes the factor \(c^3\)
from the twist while retaining coprimality at its prime divisors.

For an unramified rational twist, the square part \(b^2\) forces a
zero theta coefficient unless \(b=1\). Otherwise the normalized Gauss sum for this integer equals \(1\), by \eqref{eq:rational-gauss}.
The residue computation in the proof of \cite[Lemma 6.7]{DDFDS}
then contributes
\[
a^{-\frac{1}{3}}
\prod_{\substack{\varpi\mid ac\\(\varpi,3)=1}}
\frac{1-N(\varpi)^{-1}}{1-N(\varpi)^{-2}}
=a^{-\frac{1}{3}}\prod_{\substack{p\mid ac\\p\ne3}}\kappa_p.
\]
This includes primes occurring only in the removed cube.
Combining the ramified and unramified factors proves the formula.
\end{proof}

\subsection{Some results on multivariable complex functions}

We use Bochner's tube theorem in the form of \cite[Theorem 2.12]{GZ}:
a holomorphic function on a connected tube domain extends holomorphically
to its convex hull. To control its growth after continuation, we use
the comparison principle in \cite[Proposition 2.13]{GZ}. This transfers
a bound by the absolute value of a holomorphic function to the convex
hull. In our application, the comparison function is a polynomial
with no zeros on the relevant tube domain.

\section{Proof of Theorem~\ref{thm:main}}\label{sec:first-moment-proof}

\subsection{Initial treatment}

For large real parts, define
\begin{equation}
A_3(s,w)=\sum_{(q,3)=1}q^{-s}
 \sideset{}{^*}\sum_{\substack{\chi\bmod q\\\chi^3=\chi_0}}L(w,\chi).
\label{eq:A3}
\end{equation}
As in \cite{GZ}, a prime on a sum over primary \(n\) means that \(n\)
is squarefree and has no rational prime divisor. In such a sum we
include \(n=1\). Consequently
\[
A_3(s,w)+\zeta(w)
=\sideset{}{'}\sum_{n\ {\rm primary}}
\frac{L(w,\widehat\chi_{3,n})}{N(n)^s}.
\]
The explicit term \(\zeta(w)\) accounts for conductor \(1\), which is
excluded from \eqref{eq:A3}.

For \(c>1\), Mellin inversion and GRH give
\begin{equation}
\sum_{(q,3)=1}\sideset{}{^*}\sum_{\substack{\chi\bmod q\\
\chi^3=\chi_0}}L\left(\frac{1}{2},\chi\right)\Phi\left(\frac{q}{Q}\right)
=\frac1{2\pi i}\int_{(c)}
 A_3\left(s,\frac{1}{2}\right)\widehat\Phi(s)Q^s\,ds.
\label{eq:Mellin}
\end{equation}
\subsection{First region of absolute convergence of \texorpdfstring{\(A_3(s,w)\)}{A3(s,w)}}
\label{subsec:first-region}

Write \(\sigma=\Re s\) and \(\tau=\Re w\). For each positive integer
\(m\), let \(\psi_m\) be the Hecke character of order dividing \(3\)
whose value at an
ideal \((n)\) coprime to \(3m\), with \(n\) primary, is
\(\psi_m((n))=\left(\frac{m}{n}\right)_3\). We extend \(\psi_m\)
by zero at ideals not coprime to \(3m\); thus its Hecke
\(L\)-function has no Euler factors at primes dividing \(3m\).
Initially for \(\sigma,\tau>1\), expanding the Dirichlet
\(L\)-functions and interchanging the two absolutely convergent sums
gives
\[
A_3(s,w)+\zeta(w)=\sum_{m=1}^{\infty}\frac{F_m(s)}{m^w},
\qquad
F_m(s)=\sideset{}{'}\sum_{n\ {\rm primary}}
\frac{\psi_m((n))}{N(n)^s}.
\]

Following \cite[Section 3.2]{GZ}, we remove the restriction on
rational prime divisors by M\"obius inversion. The unrestricted
squarefree sum has Euler product
\[
\prod_{\varpi\nmid3m}
\left(1+\frac{\psi_m(\varpi)}{N(\varpi)^s}\right)
=\frac{L(s,\psi_m)}{L(2s,\psi_m^2)}.
\]
Here \(\varpi\) runs over primary primes in \(K\), and we write
\(\psi_m(\varpi)=\psi_m((\varpi))\). Since
\(\psi_m((d))=1\) for rational \(d\) coprime to \(3m\), the
inversion gives
\[
\begin{aligned}
F_m(s)
&=\frac{L(s,\psi_m)}{L(2s,\psi_m^2)}
  \sum_{\substack{d\ge1\\(d,3m)=1}}
  \frac{\mu_{\Q}(d)}{d^{2s}}
  \prod_{\varpi\mid d}
  \left(1+\frac{\psi_m(\varpi)}{N(\varpi)^s}\right)^{-1}\\
&=\frac{L(s,\psi_m)}{L(2s,\psi_m^2)}
  \prod_{p\nmid3m}
  \left(1-\frac{1}{p^{2s}}
  \prod_{\varpi\mid p}
  \left(1+\frac{\psi_m(\varpi)}{N(\varpi)^s}\right)^{-1}\right).
\end{aligned}
\]
In the last product, \(p\) runs over rational primes. On every half-plane
\(\sigma\ge\frac12+\delta\), the local factors of the last product
are \(1+O_\delta(p^{-2\sigma})\), uniformly in \(m\). Hence that
product converges normally, and the reciprocal of
\(L(2s,\psi_m^2)\) is bounded there by its absolutely convergent
Euler product. The only possible pole of the numerator is at
\(s=1\), when \(\psi_m\) is principal.

The Lindel\"of consequence of GRH for the numerator, together with
the bounds for the omitted Euler factors, gives, for every
\(\delta,\eta>0\),
\[
(s-1)F_m(s)\ll_{\delta,\eta}
m^{\eta}(1+|\Im s|)^A,
\qquad \sigma\ge\frac12+\delta,
\]
with \(\sigma\) restricted to a fixed bounded interval. Choosing
\(\eta<\tau-1\), locally uniformly in \(\tau>1\), makes the
continued \(m\)-series normally convergent after multiplication by
\(s-1\), since
\[
\sum_{m=1}^{\infty}\frac{|(s-1)F_m(s)|}{m^{\tau}}
\ll_{\delta,\eta}(1+|\Im s|)^A
\sum_{m=1}^{\infty}m^{-\tau+\eta}<\infty.
\]
It follows that \(A_3(s,w)\) continues meromorphically to
\begin{equation}
S_1=\left\{(s,w):\sigma>\frac{1}{2},\ \tau>1\right\}.
\label{eq:S1}
\end{equation}
Since \(\zeta(w)\) is holomorphic and bounded on smaller tubes in
\(\tau>1\), the function \((s-1)A_3(s,w)\) is holomorphic in
\(S_1\) and has polynomial growth in the imaginary parts there,
locally uniformly in the real parts.

\subsection{Second region of absolute convergence of \texorpdfstring{\(A_3(s,w)\)}{A3(s,w)}}
\label{subsec:second-region}

We first use the defining series \eqref{eq:A3}. On a fixed bounded
interval of \(\tau=\Re w\), GRH and the primitive functional
equation give
\[
L(w,\chi)\ll_{\eps}
q^{\max\{0,\frac{1}{2}-\tau\}+\eps}(1+|\Im w|)^A.
\]
Here GRH gives the bound when \(\tau\ge\frac{1}{2}\), while the
functional equation reduces the case \(\tau<\frac{1}{2}\) to that
range. There are \(O_{\eps}(q^{\eps})\) primitive cubic characters
of conductor \(q\). Hence
\[
\sum_{(q,3)=1}q^{-\sigma}
\sideset{}{^*}\sum_{\substack{\chi\bmod q\\\chi^3=\chi_0}}
|L(w,\chi)|
\ll_{\eps}(1+|\Im w|)^A
\sum_{q\ge1}q^{-\sigma+\max\{0,\frac12-\tau\}+2\eps}.
\]
Choosing \(\eps\) sufficiently small gives convergence whenever
\(\sigma>1+\max\{0,\frac12-\tau\}\), or equivalently in
\begin{equation}
S_2=\left\{(s,w):\sigma>1,\ \sigma+\tau>\frac{3}{2}\right\}.
\label{eq:S2}
\end{equation}
The convergence is locally uniform, so \(A_3(s,w)\) is holomorphic
and polynomially bounded in the imaginary parts on smaller tubes
in \(S_2\). The exclusion of the conductor \(1\) character ensures
holomorphy at \(w=1\).

To obtain a further region, we apply the functional equation before
summing over the conductors. All cubic Dirichlet characters are
even, and \cite[Lemma 2.6]{GZ} gives
\[
\tau(\widehat\chi_{3,n})
=\overline{\left(\frac{c_K}{n}\right)_3}g_K(\chi_{3,n})
=g_K(c_K,\chi_{3,n}).
\]
The conductor factor in the functional equation is
\(N(n)^{-w}\) when the Gauss sum is left unnormalized. Thus
\begin{equation}
A_3(s,w)+\zeta(w)
=\pi^{w-\frac{1}{2}}
\frac{\Gamma\left(\frac{1-w}{2}\right)}{\Gamma\left(\frac{w}{2}\right)}
 C_{3,1}(s+w,1-w),                               \label{eq:FE}
\end{equation}
where, in the notation of \cite{GZ},
\begin{equation}
C_{3,1}(s,w)=
\sum_{m=1}^\infty\frac1{m^w}
\sideset{}{'}\sum_{n\ {\rm primary}}
\overline{\left(\frac{mc_K}{n}\right)_3}
\frac{g_K(\chi_{3,n})}{N(n)^s}.
\label{eq:C31}
\end{equation}
The term \(\zeta(w)\) accounts for \(n=1\), which is included in
the primed sum. Since the cubic symbol is trivial at \(-1\),
\(C_{3,-1}=C_{3,1}\).

M\"obius inversion expresses the restriction on rational prime
divisors in \eqref{eq:C31} as a sum over \(d\). Cubic reciprocity
and \eqref{eq:rational-gauss} then give
\begin{equation}
C_{3,1}(s,w)=
\sum_{\substack{m,d\ge1\\(m,d)=1,\ (d,3)=1}}
\frac{\mu_{\Q}(d)d}{m^w d^{2s}}\,
 H_{3,1}(s,mdc_K),
\qquad \Re s>\frac{3}{2},\quad \Re w>1.
\label{eq:Mobius}
\end{equation}
This is the cubic specialization of the decomposition following
\cite[(3.14)]{GZ}. The denominator \(d^{2s}\) comes from
\(N(d)=d^2\), and the numerator \(d\) is the Gauss sum
\(g_K(c_K,\chi_{3,d})\) evaluated in \eqref{eq:rational-gauss}.

We estimate the \(m\)- and \(d\)-sums jointly by localizing them
to finite scales, retaining the M\"obius coefficients inside the
average.
Choose \(U,V\in C_c^\infty(0,\infty)\), and set
\begin{equation}
W(x)=(xU(x))',\qquad
\widehat W(z)=(1-z)\widehat U(z).
\label{eq:weight}
\end{equation}
The Mellin identity follows by integration by parts; compact
support makes the boundary terms vanish. In particular,
\(\widehat W(1)=0\). The factor \(1-z\) will cancel the pole of a
principal numerator \(L(z,\chi)\)
in the joint Mellin estimate below. Its effect on the later
reconstruction is explicit, since
\(\widehat W(1-w)=w\widehat U(1-w)\).
For \(M,D>0\), define
\begin{align}
B(M,D;s)={}&\left(s-\frac{4}{3}\right)\zeta_{c_K}(3s-2)                 \notag\\
&\times
\sum_{\substack{m,d\ge1\\(m,d)=1,\ (d,3)=1}}
\mu_{\Q}(d)d\,H_{3,1}(s,mdc_K)W\left(\frac{m}{M}\right)V\left(\frac{d}{D}\right).
\label{eq:block}
\end{align}
For fixed \(M,D\), the sum is finite. The factor in front is
precisely the one in Lemma~\ref{lem:Hbound}, so \(B(M,D;s)\) is
holomorphic in \(\Re s>\frac{2}{3}\). Bounds on two vertical lines
can therefore be interpolated before the scale integrations in
\eqref{eq:reconstruction} recover the double series. The weight
\(W\) is auxiliary and does not change the external weight
\(\Phi\) in Theorem~\ref{thm:main}.

\begin{lemma}\label{lem:blocks}
Assume GRH. For fixed \(0<\kappa<\frac{1}{6}\), on
\(\Re s=\frac{3}{2}+\kappa\), and for every \(\delta>0\),
\begin{equation}
B(M,D;s)\ll_{U,V,\kappa,\delta}
M^{\frac{1}{2}+\delta}D^{\frac{3}{2}+\delta}(1+|\Im s|)^A.
\label{eq:right-block}
\end{equation}
On \(\Re s=\frac{2}{3}+\kappa\), for every \(\eps>0\),
\begin{equation}
B(M,D;s)\ll_{U,V,\kappa,\eps}
M^{\frac{5}{2}-\Re s+\eps}D^{\frac{7}{2}-\Re s+\eps}
(1+|\Im s|)^A.
\label{eq:left-block}
\end{equation}
Consequently, on compact substrips of \(\frac{2}{3}<\Re s<\frac{3}{2}\),
\begin{equation}
B(M,D;s)\ll
M^{\frac{29}{10}-\frac{8\Re s}{5}+\eps}
D^{\frac{39}{10}-\frac{8\Re s}{5}+\eps}
(1+|\Im s|)^A.
\label{eq:block-interpolation}
\end{equation}
These estimates hold for \(M,D\ge1\). Values of \(M\) or \(D\)
below \(1\) do not affect the convergence arguments below.
\end{lemma}

\begin{proof}
For \(\Re s>\frac{3}{2}\), the defining series of \(H_{3,1}\) gives
\begin{align*}
H_{3,1}(s,mdc_K)
&=\sum_{\substack{n\ {\rm primary}\\(n,md)=1}}
\frac{g_K(mdc_K,\chi_{3,n})}{N(n)^s}\\
&=\sum_{\substack{n\ {\rm primary}\\n\ {\rm squarefree}}}
\frac{g_K(c_K,\chi_{3,n})}{N(n)^s}
\overline{\widehat\chi}_{3,n}(md).
\end{align*}
Indeed, coprimality forces the Gauss sum to vanish unless \(n\) is
squarefree, and changing the additive twist contributes
\(\overline{\widehat\chi}_{3,n}(md)\). Extending this character by
zero absorbs the condition \((n,md)=1\). Substitution in
\eqref{eq:Mobius} and absolute convergence yield, for
\(\Re s>\frac{3}{2}\) and \(\Re w>1\),
\begin{equation}
C_{3,1}(s,w)=
\sum_{\substack{n\ {\rm primary}\\n\ {\rm squarefree}}}
\frac{g_K(c_K,\chi_{3,n})}{N(n)^s}
\sum_{\substack{m,d\ge1\\(m,d)=1,\ (d,3)=1}}
\frac{\mu_{\Q}(d)\overline{\widehat\chi}_{3,n}(md)}
 {m^w d^{2s-1}}.
\label{eq:C31-open}
\end{equation}
The outer sum is over all squarefree primary \(n\), including those
with rational prime divisors; the original restriction is encoded
by the M\"obius sum over \(d\). Similarly, \eqref{eq:block} becomes
\begin{align}
B(M,D;s)
={}&\left(s-\frac{4}{3}\right)\zeta_{c_K}(3s-2)
\sum_{\substack{n\ {\rm primary}\\n\ {\rm squarefree}}}
\frac{g_K(c_K,\chi_{3,n})}{N(n)^s}                  \notag\\
&\quad\times
\sum_{\substack{m,d\ge1\\(m,d)=1,\ (d,3)=1}}
\mu_{\Q}(d)d\,\overline{\widehat\chi}_{3,n}(md)
W\left(\frac{m}{M}\right)V\left(\frac{d}{D}\right).
\label{eq:block-open}
\end{align}
We first estimate the inner sum uniformly for the possibly
imprimitive Dirichlet character
\(\chi=\overline{\widehat\chi}_{3,n}\), and then sum over \(n\).

For \(\Re z,\Re v>1\), Euler multiplication gives
\begin{align}
\sum_{\substack{m,d\ge1\\(m,d)=1,\ (d,3)=1}}
\frac{\mu_{\Q}(d)\chi(md)}{m^z d^v}
={}&
\frac{L(z,\chi)}
 {L(v,\chi)(1-\chi(3)3^{-v})}                       \notag\\
&\times
\prod_{p\ne3}
\left(1+\frac{\chi(p)^2p^{-z-v}}
 {1-\chi(p)p^{-v}}\right).
\label{eq:ratio}
\end{align}
At a prime \(p\ne3\), the possibilities \(p\nmid d\) and
\(p\mid d\) give the local factor
\[
\frac1{1-\chi(p)p^{-z}}-\chi(p)p^{-v}.
\]
At \(3\), only the first possibility occurs, giving
\((1-\chi(3)3^{-z})^{-1}\). This accounts for the separate factor
at \(3\) in \eqref{eq:ratio}. The correction product converges
normally when \(\Re v>0\) and \(\Re(z+v)>1\).

To apply this identity to the smoothed sum, write
\(dV(d/D)=D(d/D)V(d/D)\). Since the Mellin transform of
\(xV(x)\) at \(v\) is \(\widehat V(v+1)\), double Mellin
inversion gives, for \(c_z,c_v>1\),
\begin{align*}
&\sum_{\substack{m,d\ge1\\(m,d)=1,\ (d,3)=1}}
\mu_{\Q}(d)d\,\chi(md)
W\left(\frac{m}{M}\right)V\left(\frac{d}{D}\right)\\
&\quad=
\frac{D}{(2\pi i)^2}
\int_{(c_z)}\int_{(c_v)}
\widehat W(z)\widehat V(v+1)M^zD^v
\frac{L(z,\chi)}{L(v,\chi)(1-\chi(3)3^{-v})}\\
&\qquad\qquad\times
\prod_{p\ne3}\left(1+
\frac{\chi(p)^2p^{-z-v}}{1-\chi(p)p^{-v}}\right)
\,dv\,dz.
\end{align*}
Shift first the \(z\)-line and then the \(v\)-line to
\(\frac{1}{2}+\delta\). The correction product remains normally
convergent throughout these shifts. Under GRH,
\(L(v,\chi)^{-1}\) has no poles in the region traversed; for a
principal character, its value at \(v=1\) is a zero.
The possible numerator pole at \(z=1\) is canceled by
\(\widehat W(z)=(1-z)\widehat U(z)\).
Thus no residues occur in the shifts, including when \(\chi\) is
principal. On the new lines the scale factors have absolute value
\(M^{\frac{1}{2}+\delta}D^{\frac{3}{2}+\delta}\); the extra power
of \(D\) comes from the weight \(d\). The GRH estimates and the
rapid decay of the Mellin transforms now give
\begin{equation}
\sum_{\substack{m,d\ge1\\(m,d)=1,\ (d,3)=1}}
\mu_{\Q}(d)d\,\overline{\widehat\chi}_{3,n}(md)
 W\left(\frac{m}{M}\right)V\left(\frac{d}{D}\right)
\ll_{\delta,\eta,U,V}
M^{\frac{1}{2}+\delta}D^{\frac{3}{2}+\delta}N(n)^\eta
\label{eq:inner-average}
\end{equation}
for every \(\eta>0\). Here the modulus of \(\chi\) is at most
\(N(n)\), and the imprimitive Euler factors are absorbed in
\(N(n)^\eta\). The conductor loss \(\eta\) can be chosen
independently of the contour displacement \(\delta\).

On \(\Re s=\frac{3}{2}+\kappa\), choose \(\eta<\kappa\).
Since \(|g_K(c_K,\chi_{3,n})|=N(n)^{\frac{1}{2}}\) for
squarefree \(n\), substitution into \eqref{eq:block-open} bounds
the outer sum by
\[
\sum_{\substack{n\ {\rm primary}\\n\ {\rm squarefree}}}
\frac{|g_K(c_K,\chi_{3,n})|N(n)^\eta}
 {N(n)^{\frac{3}{2}+\kappa}}
\le\sum_{n\ {\rm primary}}N(n)^{-1-\kappa+\eta}
\ll\zeta_K(1+\kappa-\eta)<\infty.
\]
The factor \(\left(s-\frac{4}{3}\right)\zeta_{c_K}(3s-2)\) is
polynomially bounded on this line. This proves
\eqref{eq:right-block}.

For the left boundary, Lemma~\ref{lem:Hbound} applies directly to
each term in \eqref{eq:block}. Compact support restricts the sum
to \(m\asymp M\) and \(d\asymp D\), with \(O(MD)\) pairs.
Thus, on \(\sigma=\frac{2}{3}+\kappa\),
\begin{align*}
|B(M,D;s)|
&\ll (1+|\Im s|)^A
\sum_{\substack{m\asymp M\\d\asymp D}}
d(md)^{\frac{3}{2}-\sigma+\eps}\\
&\ll M^{\frac{5}{2}-\sigma+\eps}
D^{\frac{7}{2}-\sigma+\eps}(1+|\Im s|)^A,
\end{align*}
which is \eqref{eq:left-block}.

We interpolate between the actual boundary lines
\[
a=\frac{2}{3}+\kappa,\qquad b=\frac{3}{2}+\kappa.
\]
With small losses \(\ell,\delta>0\), define the affine function
\begin{equation}
E(\sigma)=
\frac{b-\sigma}{b-a}\left(\frac{5}{2}-a+\ell\right)
+\frac{\sigma-a}{b-a}\left(\frac{1}{2}+\delta\right).
\label{eq:affine-exponent}
\end{equation}
The boundary exponents of \(M\) are \(E(a)\) and \(E(b)\),
and those of \(D\) are larger by \(1\).
For fixed \(M,D\), the holomorphic function \(B(M,D;s)\) has
polynomial growth throughout the closed strip. This follows from
\eqref{eq:completed-input} and the finite coprimality formula,
using an upper parameter slightly larger than \(b\) in
\eqref{eq:completed-input}.
Extend \(E\) affinely to complex arguments. Dividing
\(B(M,D;s)\) by
\(M^{E(s)}D^{E(s)+1}(s+B_0)^N\), with \(B_0>0\) and \(N\)
sufficiently large, we may apply Hadamard's three lines theorem.
It follows that
\[
B(M,D;s)\ll
M^{E(\sigma)}D^{E(\sigma)+1}(1+|\Im s|)^A,
\qquad a\le\sigma\le b.
\]
As \(\kappa,\ell,\delta\to0^+\),
\[
E(\sigma)\longrightarrow
\frac{\frac32-\sigma}{\frac56}\,\frac{11}{6}
+\frac{\sigma-\frac23}{\frac56}\,\frac12
=\frac{29}{10}-\frac{8\sigma}{5}
\]
uniformly on bounded intervals. On a fixed compact substrip of
\(\frac{2}{3}<\Re s<\frac{3}{2}\), we may therefore choose
\(\kappa,\ell,\delta\) sufficiently small that \(a<\sigma<b\) and
\(E(\sigma)\le\frac{29}{10}-\frac{8\sigma}{5}+\eps\).
Since \(M,D\ge1\), the preceding bound proves
\eqref{eq:block-interpolation}.
All constants can be controlled by finitely many derivative
seminorms of \(U,V\), as required for the uniform vertical bounds
below.
\end{proof}

We now determine the region supplied by Lemma~\ref{lem:blocks}.
To recover \(C_{3,1}(s+w,1-w)\), the scale integrals must restore
the factors \(m^{w-1}\) and \(d^{-2(s+w)}\) in
\eqref{eq:Mobius}. We therefore integrate \(B(M,D;s+w)\) against
\(M^{w-1}D^{-2(s+w)}\frac{dM}{M}\frac{dD}{D}\).
The exact reconstruction identity is established in
\eqref{eq:reconstruction} below. For now, applying
\eqref{eq:block-interpolation} with real part \(\sigma+\tau\)
gives the following total scale exponents, apart from arbitrarily
small positive losses:
\begin{equation}
\begin{aligned}
\left(\frac{29}{10}-\frac85(\sigma+\tau)\right)+\tau-1
&=\frac{19}{10}-\frac85\sigma-\frac35\tau,\\
\left(\frac{39}{10}-\frac85(\sigma+\tau)\right)-2(\sigma+\tau)
&=\frac{39}{10}-\frac{18}{5}(\sigma+\tau).
\end{aligned}
\label{eq:scale-exponents}
\end{equation}
Since \(\int_1^\infty x^\alpha\frac{dx}{x}<\infty\) exactly when
\(\alpha<0\), convergence at infinity requires both exponents
to be negative. These two conditions are
\[
16\sigma+6\tau>19,\qquad \sigma+\tau>\frac{13}{12}.
\]
The interpolation estimate also requires
\(\frac23<\sigma+\tau<\frac32\); its lower bound is already
implied by the second inequality above. These conditions give
\begin{equation}
S_{3,3}=\left\{(s,w):16\sigma+6\tau>19,\
\frac{13}{12}<\sigma+\tau<\frac{3}{2}\right\}.
\label{eq:S3}
\end{equation}
We combine this with \(S_1\) and \(S_2\) to obtain the following
continuation and bounds.

\begin{proposition}\label{prop:continuation}
Assume GRH. The function
\[
(s-1)\left(s+w-\frac{4}{3}\right)A_3(s,w)
\]
extends holomorphically to
\begin{equation}
S_{4,3}=\left\{(s,w):\sigma>\frac{1}{2},\ \sigma+\tau>\frac{13}{12},\
14\sigma+9\tau>16\right\}.
\label{eq:S4}
\end{equation}
It has polynomial growth in the two imaginary parts, locally
uniformly in the real parts. In particular, \(A_3\left(s,\frac{1}{2}\right)\) is
meromorphic for \(\Re s>\frac{23}{28}\), with at most simple poles at
\(s=1\) and \(s=\frac{5}{6}\), and is polynomially bounded on closed
vertical strips away from those poles.
\end{proposition}

\begin{proof}
\emph{Mellin reconstruction.}
We first recover \(C_{3,1}(s+w,1-w)\) from the finite sums
\(B(M,D;s+w)\). For fixed positive integers \(m,d\), the changes of
variables \(x=m/M\) and \(y=d/D\) give
\[
\begin{aligned}
\int_0^\infty W(m/M)M^{w-1}\frac{dM}{M}
&=m^{w-1}\int_0^\infty W(x)x^{-w}\,dx
 =w\,m^{w-1}\widehat U(1-w),\\
\int_0^\infty V(d/D)D^{-2(s+w)}\frac{dD}{D}
&=d^{-2(s+w)}\int_0^\infty V(y)y^{2(s+w)-1}\,dy
 =d^{-2(s+w)}\widehat V(2s+2w).
\end{aligned}
\]
In the first identity we used \eqref{eq:weight} at \(z=1-w\).
These identities may be summed term by term when
\(\sigma+\tau>\frac32\) and \(\tau<0\). Indeed, absolute
convergence of the defining Gauss sum series gives
\(H_{3,1}(s+w,mdc_K)\ll_{\sigma+\tau}1\), uniformly in \(m,d\).
After taking absolute values, the two scale integrals are therefore
bounded by a constant times
\[
\sum_{m=1}^\infty m^{\tau-1}
\sum_{d=1}^\infty d^{1-2(\sigma+\tau)}<\infty.
\]
Fubini's theorem and \eqref{eq:Mobius} now give the exact identity
\begin{align}
&\int_0^\infty\int_0^\infty
 B(M,D;s+w)M^{w-1}D^{-2(s+w)}
 \frac{dM}{M}\frac{dD}{D}                           \notag\\
&\quad=
w\widehat U(1-w)\widehat V(2s+2w)
\left(s+w-\frac{4}{3}\right)\zeta_{c_K}(3s+3w-2)
 C_{3,1}(s+w,1-w).
\label{eq:reconstruction}
\end{align}
Denote the integral on the left by \(I_{U,V}(s,w)\).

\emph{Continuation of the scale integral.}
The exponents in \eqref{eq:scale-exponents} are strictly negative
on \(S_{3,3}\). On any fixed compact set of real parts in this
region, they are at most \(-\eta\) for some \(\eta>0\).
Choose the loss in \eqref{eq:block-interpolation} smaller than
\(\eta/2\). The integrand at large scales is then bounded by
\[
M^{-\eta/2}D^{-\eta/2}
(1+|\Im(s+w)|)^A
\]
with respect to \(\frac{dM}{M}\frac{dD}{D}\).
Compact support makes \(B\) vanish when either scale is below a
fixed positive number; a fixed rescaling of the weights extends
the same estimates to the remaining scales below \(1\).
Normal convergence therefore makes \(I_{U,V}\) holomorphic in
\(S_{3,3}\). Moreover,
\[
19<16\sigma+6\tau=16(\sigma+\tau)-10\tau<24-10\tau,
\]
so \(\tau<\frac12\) throughout this region.

To connect this integral with the original identity, we keep
\(\kappa,\ell,\delta>0\) in \eqref{eq:affine-exponent}.
The corresponding scale exponents are
\[
E(\sigma+\tau)+\tau-1,
\qquad E(\sigma+\tau)+1-2(\sigma+\tau).
\]
At
\[
\sigma+\tau=\frac32+\frac\kappa2,
\qquad \tau=-\frac14,
\]
write
\[
E_*=
\frac{3\kappa}{5}\left(\frac{11}{6}-\kappa+\ell\right)
+\left(1-\frac{3\kappa}{5}\right)
 \left(\frac12+\delta\right).
\]
The two exponents are \(E_*-\frac54\) and \(E_*-2-\kappa\),
which tend to \(-\frac34\) and \(-\frac32\), respectively, as
\(\kappa,\ell,\delta\to0^+\). They are therefore strictly negative
for sufficiently small parameters. The integral converges in an
open neighborhood of this point, where
\(\sigma=\frac74+\frac\kappa2>1\); this neighborhood lies in
\(S_2\) and in the initial domain of \eqref{eq:reconstruction}.

For fixed small parameters, negativity of these two affine
exponents, together with
\[
a<\sigma+\tau<b,\qquad
\sigma+\tau>\frac{13}{12},\qquad \tau<\frac12,
\]
defines a convex region. For any compact subset of \(S_{3,3}\),
small enough parameters make this region contain both that subset
and the initial open overlap. Varying the parameters thus gives
a connected tube containing \(S_{3,3}\). The resulting functions
agree on overlaps, where they are represented by the same
absolutely convergent integral.

\emph{Removal of the auxiliary Mellin factors.}
Zeros of \(\widehat U(1-w)\widehat V(2s+2w)\) do not obstruct
continuation. For two choices \((U_1,V_1)\), \((U_2,V_2)\),
\eqref{eq:reconstruction} gives on the initial overlap
\[
I_{U_1,V_1}\widehat U_2(1-w)\widehat V_2(2s+2w)
=I_{U_2,V_2}\widehat U_1(1-w)\widehat V_1(2s+2w).
\]
The identity theorem extends this equality to each connected
region above. At any point, choose \(U,V\) with both transforms
nonzero, using the modulation below. Division by their product
then defines compatible local holomorphic quotients.

Combining \eqref{eq:reconstruction} with \eqref{eq:FE} gives,
initially on the original overlap,
\[
\begin{aligned}
&\left(s+w-\frac43\right)\zeta_{c_K}(3s+3w-2)
 \bigl(A_3(s,w)+\zeta(w)\bigr)\\
&\qquad=
\frac{\pi^{w-\frac12}\Gamma\left(\frac{1-w}{2}\right)}
 {w\Gamma\left(\frac w2\right)}
\frac{I_{U,V}(s,w)}{\widehat U(1-w)\widehat V(2s+2w)}.
\end{aligned}
\]
The right-hand side is holomorphic throughout the new region.
In particular, its apparent singularity at \(w=0\) is removable,
since
\[
\frac{\pi^{w-\frac12}\Gamma\left(\frac{1-w}{2}\right)}
 {w\Gamma\left(\frac w2\right)}\longrightarrow\frac12.
\]
The numerator gamma function has no poles for \(\tau<\frac12\).
Since \(\zeta(w)\) is also holomorphic there, subtraction of the
conductor \(1\) contribution and multiplication by \(s-1\) give a
holomorphic continuation of
\begin{equation}
(s-1)\left(s+w-\frac{4}{3}\right)\zeta_{c_K}(3s+3w-2)A_3(s,w).
\label{eq:cleared}
\end{equation}
By Sections~\ref{subsec:first-region} and~\ref{subsec:second-region},
this function is already holomorphic in \(S_1\) and \(S_2\).
Both of those regions satisfy \(\sigma+\tau>\frac32\), so the
additional zeta factor introduces no pole. Their open intersection,
together with the connecting regions above, gives a connected tube
on which \eqref{eq:cleared} is holomorphic and which contains
\(S_1\), \(S_2\), and \(S_{3,3}\).

\emph{Vertical bounds and the convex hull.}
We first bound \eqref{eq:cleared} on this connected tube. The bounds
on \(S_1\) and \(S_2\) have already been established. For the new
regions, fix nonnegative, nonzero \(U_0,V_0\in C_c^\infty(1,2)\).
At a given point \((s,w)\), choose
\[
U(x)=x^{i\Im w}U_0(x),\qquad
V(x)=x^{-2i\Im(s+w)}V_0(x).
\]
Then
\[
\widehat U(1-w)=\widehat U_0(1-\tau)>0,
\qquad
\widehat V(2s+2w)=\widehat V_0(2\sigma+2\tau)>0.
\]
The transforms are bounded away from zero on compact sets of real
parts, while the derivative seminorms of \(U,V\) grow polynomially
in \(|\Im s|+|\Im w|\). Lemma~\ref{lem:blocks}, the scale integral,
Stirling's formula, and the vertical-strip bounds for \(\zeta(w)\)
therefore give locally uniform polynomial growth for
\eqref{eq:cleared}. This pointwise choice is used only for the
bound; during analytic continuation the weights remain fixed on
each neighborhood.

The convex hull of the real bases of \(S_1\), \(S_2\), and
\(S_{3,3}\) is the base of \(S_{4,3}\). Its two finite boundary
vertices are
\[
\left(\frac12,1\right),\qquad
\left(\frac54,-\frac16\right),
\]
joined by \(14\sigma+9\tau=16\), and its boundary rays have
directions \((0,1)\) and \((1,-1)\). The connected tube constructed
above contains these three regions, so Bochner's theorem extends
\eqref{eq:cleared} holomorphically to \(S_{4,3}\).
In particular, for \(h>0\) small, the points
\(\left(\frac12+h,1+h\right)\in S_1\) and
\(\left(\frac54+h,-\frac16+h\right)\in S_{3,3}\), with respective
weights \(\frac47-\frac{6h}{7}\) and \(\frac37+\frac{6h}{7}\),
have convex combination
\begin{equation}
\left(\frac{23}{28}+\frac{23h}{14},\ \frac12\right).
\label{eq:central-interior}
\end{equation}

To transfer the growth bounds, choose a bounded connected real
open set whose compact closure lies in the connected base just
constructed and whose convex hull contains a prescribed compact
subset of the base of \(S_{4,3}\). Such a set can be obtained from
finitely many neighborhoods and paths through the overlaps.
The polynomial bounds are uniform over its closure. For
sufficiently large \(B_0,N\) and \(C>0\),
\[
C(B_0+s)^N(B_0+w)^N
\]
is nonzero on the tube over its convex hull and dominates the
absolute value of \eqref{eq:cleared} on the original tube.
The comparison principle \cite[Proposition 2.13]{GZ}, or
\cite[Proposition C.5]{Cech}, transfers this domination to the
convex hull. This proves polynomial growth on \(S_{4,3}\),
locally uniformly in the real parts.

Finally, throughout \(S_{4,3}\),
\[
\Re(3s+3w-2)>\frac54>1.
\]
The zeta factor and its reciprocal are therefore holomorphic and
bounded on smaller tubes. Dividing \eqref{eq:cleared} by
\(\zeta_{c_K}(3s+3w-2)\) proves the asserted continuation and
growth. On the central line \(w=\frac12\), the inequalities in
\eqref{eq:S4} reduce to \(\Re s>\frac{23}{28}\), and the two
remaining possible polar factors are \(s-1\) and
\(s-\frac56\). This gives the final assertion.
\end{proof}

\subsection{Residue of \texorpdfstring{\(A_3(s,w)\) at \(s=1\)}{A3(s,w) at s=1}}

The residue calculation in \cite{GZ} gives, initially for \(\Re w>1\),
\begin{equation}
\Res_{s=1}A_3(s,w)
=\frac{3r_K\zeta(3w)}{4\zeta_K(2)}
\prod_{p\ne3}
\left(1-\kappa_pp^{-2}
 +\bigl(\kappa_p(1+p^{-2})-1\bigr)p^{-3w}\right).
\label{eq:first-residue}
\end{equation}
Indeed, before extracting \(\zeta(3w)\), the local factor at \(p\ne3\)
is
\[
1-\kappa_pp^{-2}
+\kappa_p\frac{p^{-3w}}{1-p^{-3w}},
\]
and the factor at \(3\) is
\(\left(\frac{3}{4}\right)(1-3^{-3w})^{-1}\).
This is the Euler product form of \cite[(3.19) and (3.20)]{GZ}.
It continues to the range needed here; at \(w=\frac{1}{2}\) it gives
\eqref{eq:C3}.

\subsection{Completion of the proof}

We now shift the line of integration in \eqref{eq:Mellin}.
Choose \(0<\eps<\frac{1}{84}\). Proposition~\ref{prop:continuation}
permits the line to be moved to \(\Re s=\frac{23}{28}+\eps\).
Its polynomial growth bound, uniform across the intervening
closed strip after removing the poles, combines with
\[
\widehat\Phi(\sigma+it)\ll_{\Phi,J,A}(1+|t|)^{-A}
\quad(\sigma\in J)
\]
for every compact real interval \(J\) and every \(A>0\).
The horizontal integrals tend to zero, and the new vertical
integral is \(O_{\Phi,\eps}\left(Q^{\frac{23}{28}+\eps}\right)\).
The residue at \(s=1\) is given by \eqref{eq:first-residue} at
\(w=\frac{1}{2}\). The residue at \(s=\frac{5}{6}\) is evaluated
in Section~\ref{sec:first-secondary-residue}; equation~\eqref{eq:central-residue}
identifies it with \(C_3'\), and \eqref{eq:C3prime} gives its
explicit value and strict negativity. These two residues give \eqref{eq:main}.
Larger values of \(\eps\) follow from the result with a smaller
positive value. This completes the proof of
Theorem~\ref{thm:main}.

\section{Evaluation and nonvanishing of the secondary terms}
\label{sec:secondary-terms}

We collect here the arithmetic calculations of the secondary terms
in Theorems~\ref{thm:main} and~\ref{thm:second-moment}.
The first moment gives a simple pole with residue \(C_3'<0\);
the second gives a double pole with leading Laurent coefficient
\(A_{\mathrm{sf}}>0\). We begin by fixing the correspondence between
the notation used for these two families and that of \cite{deFaveri}.

\subsection{Notation and comparison with de Faveri}
\label{sec:notation-comparison}

The two secondary terms arise from different families.  The first
moment in \eqref{eq:A3} averages primitive cubic Dirichlet characters
of rational conductor \(q\), using the weight \(\Phi(q/Q)\).
Its parametrization by primary elements imposes squarefreeness and
the absence of rational prime divisors.  The second moment in
\eqref{eq:second-moment-family} averages all squarefree ideals
\(\mathfrak a\) coprime to \(3\), using \(W(N\mathfrak a/X)\).

The latter is precisely the specialization of
\cite[Theorem~1.3 and Example~2.2]{deFaveri} with
\[
 n=3,\qquad K=\Q(\zeta_3)=\Q(\omega),\qquad
 \mathfrak c=(9)=(\lambda^4),\qquad
 S_{\mathrm{dF}}=\{(\lambda)\},\qquad \lambda=c_K=\sqrt{-3}.
\]
Our convention \(S=S_\infty\cup S_{\mathrm{dF}}\) gives the same
\(I(S)\) and the same omitted finite Euler factors.
The prime \((\lambda)\) is the only ramified prime of \(K/\Q\),
so the additional hypothesis on \(S\) in that theorem is satisfied.
Writing \(\mathfrak a=(a)\), \(\mathfrak b=(b)\) with primary
generators, Example~2.2 gives
\begin{equation}\label{eq:notation-character-dictionary}
 \chi^{\mathrm{dF}}_{\mathfrak a}(\mathfrak b)
 =\left(\frac{a}{b}\right)_3
 =\left(\frac{b}{a}\right)_3
 =\chi_{3,a}(b)
 \qquad ((ab,3)=1,\ (a,b)=1).
\end{equation}
The superscript \(\mathrm{dF}\) is used here only to distinguish
the conventions: \cite[Section~2.5.1]{deFaveri} writes
\(\widetilde\chi_{\mathfrak a}\) for the primitive character
inducing \(\chi^{\mathrm{dF}}_{\mathfrak a}\).  This primitive
character is our \(\chi_a\) in functional equations.  For squarefree
\(\mathfrak a\), their partial functions agree exactly:
\[
 L^{S_{\mathrm{dF}}}(s,\chi^{\mathrm{dF}}_{\mathfrak a})
 =L^S(s,\chi_a)
 =L(s,\chi_a)\bigl(1-\chi_a((\lambda))3^{-s}\bigr),
\]
where the last factor is \(1\) if \(\chi_a\) is ramified at
\(\lambda\).  Thus the second moment retains the partial
\(L\)-function used in de Faveri's theorem.
On rational integers \(m\) coprime to \(3a\), the original Dirichlet
character satisfies \(\widehat\chi_{3,a}(m)=\chi_{3,a}(m)\).
This equality of residue-symbol values does not identify the Dirichlet
\(L\)-function, summed over integers, with the Hecke \(L\)-function,
summed over ideals of \(K\).
His \(\mu(\mathfrak a)\) is our \(\mu_K(\mathfrak a)\), and
his \(\mathfrak a^*\) is the squarefree kernel of the cubefree
part, denoted \(a_1\) below.  In particular,
\[
 \operatorname{Cond}(\chi_a)=\mathfrak c_a\mathfrak a^*,
 \qquad \mathfrak c_a\mid\mathfrak c,
 \qquad \mathfrak a^*=\mathfrak a\quad
       \text{if }\mathfrak a\text{ is squarefree}.
\]
The scale in Theorem~1.3 of \cite{deFaveri} is \(N\mathfrak a/X\),
with the same \(X\) as here; replacing \(N\mathfrak a\) by the
conductor norm would define a different weighted moment.

The roles of the analytic variables are as follows:
\[
\begin{array}{c|c|c}
 & A_3(s,w) & Z_1^S(s_1,s_2,w)\\ \hline
 \text{size Mellin variable} & s & w\\
 \text{arguments of the }L\text{-functions} & w & s_1,s_2\\
 \text{central slice} & A_3(s,\frac12)
              & Z_1^S(\frac12,\frac12,w)\\
 \text{secondary pole} & s=\frac56 & w=\frac56
\end{array}
\]
The notation on the right agrees with the central slice used in
\cite[Section~9]{deFaveri}, following \cite{DIPP}.
This is a correspondence of variable roles, not an identity between
the two series. In the second-moment residue calculation we put
\(u=1-w\), so the secondary point \(w_0=\frac56\) corresponds to
\(u_0=\frac16\). Similarly, \(q\) in \eqref{eq:A3} is a rational
conductor, whereas \(q=Np\) in the second moment's local calculations
is the norm of a prime ideal \(p\).  Rational primes in the first
moment and prime ideals in the second moment are summed separately.

The second moment's nonnegative weight \(W\) is also distinct from
the auxiliary kernel in \eqref{eq:weight}.  That earlier kernel is
\((xU(x))'\) and has Mellin transform zero at \(1\); the present
weight is nonzero and satisfies \(\widehat W(\frac56)>0\).
The support condition \(\operatorname{supp}W\subset(1,2)\) in
\cite[Theorem~1.3]{deFaveri} extends to our weight class by a finite
nonnegative smooth partition and rescaling: for a piece \(W_j\)
supported in \((c_j,2c_j)\), apply that theorem to
\(V_j(t)=W_j(c_jt)\) at scale \(c_jX\).

Both calculations use the same Patterson normalization.  In the
notation of Lemma~\ref{lem:Hresidue} and \eqref{cub:rho-constant},
\begin{equation}\label{eq:notation-residue-dictionary}
 \begin{gathered}
 g_3(r,c)=g_K(c_Kr,\chi_{3,c}),\qquad
 \rho_0=27c_0
 =\frac{(2\pi)^{\frac53}}
 {8\,3^{\frac32}\Gamma(\frac23)\zeta_K(2)},\\
 \zeta_\lambda=\zeta_{c_K}=\zeta_K^S,\qquad
 \kappa_\lambda=\kappa_S
 =\Res_{s=1}\zeta_K^S(s)=\frac23r_K.
 \end{gathered}
\end{equation}
The global zeta residue \(\kappa_\lambda\) differs from the local
rational-prime weight \(\kappa_p\) in \eqref{eq:kappap}.
For a positive rational integer \(r\), let \(h\) be the squarefree
product of the prime ideals outside \(S\) dividing \((r)\).
The coprimality restriction in \(H_{3,1}\) then gives the precise
bridge to the raw series \(\psi\) and its restricted form \(\psi_h\):
\begin{equation}\label{eq:notation-restricted-residue}
 \begin{split}
 H_{3,1}(v,rc_K)&=\psi_h(r,v),\\
 \Res_{v=\frac43}H_{3,1}(v,rc_K)
 &=\beta_h\rho(r),\qquad
 \beta_h=\prod_{p\mid h}(1+(Np)^{-1})^{-1}.
 \end{split}
\end{equation}
The second equality is proved in Lemma~\ref{cub:modulus-removal-lemma}.
It identifies the common factor \(\rho_0\) in
\eqref{eq:Hresidue} and \eqref{eq:secondary-residue}, while keeping
track of the different coprimality restrictions.

Finally, the coefficient to be determined in de Faveri's notation is
the leading coefficient of \(Q_W\).  If
\[
 Z_1^S\left(\frac12,\frac12,w\right)
 =\frac{A_{\mathrm{sf}}}{(w-\frac56)^2}
  +\frac{B_{\mathrm{sf}}}{w-\frac56}+O(1),
\]
Mellin inversion gives the exact dictionary
\begin{equation}\label{eq:notation-secondary-polynomial}
 Q_W(t)=A_{\mathrm{sf}}\widehat W\left(\frac56\right)t
       +A_{\mathrm{sf}}\widehat W'\left(\frac56\right)
       +B_{\mathrm{sf}}\widehat W\left(\frac56\right).
\end{equation}
Thus \(b_W=A_{\mathrm{sf}}\widehat W'(\frac56)
+B_{\mathrm{sf}}\widehat W(\frac56)\) in
Theorem~\ref{thm:second-moment}.  We prove
\(A_{\mathrm{sf}}>0\), so \(Q_W\) has degree exactly one.
The first moment instead has the simple-pole coefficient
\(C'_3=\Res_{s=\frac56}A_3(s,\frac12)<0\).
The signs therefore describe different Laurent coefficients in
different moment problems.  The first moment asymptotic assumes
GRH; the second moment result is unconditional.

\subsection{The secondary residue in the first moment}
\label{sec:first-secondary-residue}

We now compute the residue that gives the term of order
\(Q^{\frac{5}{6}}\). We first evaluate the residue of \(C_{3,1}(s,w)\)
at \(s=\frac{4}{3}\) in a half-plane of absolute convergence, and then
continue the resulting expression to \(w=\frac{1}{2}\).

For \(\Re w>\frac{7}{6}\), Lemma~\ref{lem:Hbound} permits us to take
the residue in \eqref{eq:Mobius} term by term, giving
\begin{equation}
\Res_{s=\frac{4}{3}}C_{3,1}(s,w)
=\sum_{\substack{m,d\ge1\\(m,d)=1,\ (d,3)=1}}
 \frac{\mu_{\Q}(d)}{m^w d^{\frac{5}{3}}}\,
 \Res_{s=\frac{4}{3}}H_{3,1}(s,mdc_K).
\label{eq:residue-sum}
\end{equation}
Indeed, on a sufficiently small circle
\(\left|s-\frac{4}{3}\right|=\delta\), the summand in
\eqref{eq:Mobius} is bounded by
\[
\left|
\frac{\mu_{\Q}(d)d}{m^w d^{2s}}H_{3,1}(s,mdc_K)
\right|
\ll_{\delta,\eps}
m^{-\Re w+\frac{1}{6}+\delta+\eps}
d^{-\frac{3}{2}+3\delta+\eps}.
\]
For \(w\) in a compact subset of \(\Re w>\frac{7}{6}\), we may choose
\(\delta\) and \(\eps\) so that both exponents are less than \(-1\).
This gives locally uniform absolute convergence on the circle and
justifies the exchange of summation and contour integration.

We use Lemma~\ref{lem:Hresidue} to write
\eqref{eq:residue-sum} as an Euler product. After removing the
common constant in \eqref{eq:Hresidue}, fix a prime \(p\ne3\).
If \(p\nmid d\), the nonzero contributions with \(p\mid m\) come
from exponents \(3k\), \(k\ge1\), or \(3k+1\), \(k\ge0\).
Their respective contributions are
\(\kappa_pp^{-3kw}\) and
\(\kappa_pp^{-(3k+1)w-\frac{1}{3}}\).
If \(p\mid d\), then \((m,d)=1\) forces \(p\nmid m\), and the
corresponding contribution is
\(-p^{-\frac{5}{3}}\kappa_pp^{-\frac{1}{3}}=-\kappa_pp^{-2}\).
Together with the contribution \(1\) from \(p\nmid md\), these
terms give the local factor
\begin{equation}
\begin{aligned}
&1+\kappa_p\sum_{k=1}^{\infty}p^{-3kw}
 +\kappa_pp^{-w-\frac{1}{3}}\sum_{k=0}^{\infty}p^{-3kw}
 -\kappa_pp^{-2}\\
&\qquad=
1+\kappa_p\left(
\frac{p^{-w-\frac{1}{3}}+p^{-3w}}{1-p^{-3w}}-p^{-2}
\right).
\end{aligned}
\label{eq:local-secondary}
\end{equation}
At \(p=3\), the condition \((d,3)=1\) removes the M\"obius term,
and the ramified part of \eqref{eq:Hresidue} gives
\[
\frac{1+3^{-w-\frac{1}{3}}}{1-3^{-3w}}.
\]
These local expansions show that the residue sum is absolutely
convergent for \(\Re w>\frac{2}{3}\). To continue it to the
central point, we extract \(\zeta(3w)\) and
\(\zeta\left(w+\frac{1}{3}\right)\). The denominators
\(1-p^{-3w}\) give the first factor; multiplying the remaining
local factors by \(1-p^{-w-\frac{1}{3}}\) gives the correction
product associated with the second. At \(p=3\), the correction is
\[
\left(1-3^{-w-\frac{1}{3}}\right)
\left(1+3^{-w-\frac{1}{3}}\right)
=1-3^{-2w-\frac{2}{3}}.
\]
We obtain
\begin{align}
\Res_{s=\frac{4}{3}}C_{3,1}(s,w)
={}&\frac{(2\pi)^{\frac{5}{3}}}
 {8\,3^{\frac{3}{2}}\Gamma\left(\frac{2}{3}\right)\zeta_K(2)}
 \zeta(3w)\zeta\left(w+\frac{1}{3}\right)
 \left(1-3^{-2w-\frac{2}{3}}\right)\notag\\
&\times\prod_{p\ne3}\left(1-p^{-w-\frac{1}{3}}\right)
 \Bigl(1-\kappa_pp^{-2}+\kappa_pp^{-w-\frac{1}{3}}
 +\bigl(\kappa_p(1+p^{-2})-1\bigr)p^{-3w}\Bigr).
\label{eq:secondary-residue}
\end{align}
Since \(\kappa_p=1+O(p^{-1})\), the local factor in the remaining
product is
\[
1+O\left(
p^{-2}+p^{-\Re w-\frac{4}{3}}
+p^{-2\Re w-\frac{2}{3}}+p^{-3\Re w-1}\right),
\]
uniformly on compact subsets of \(\Re w>\frac{1}{6}\).
The product therefore converges normally in this half-plane, and
\eqref{eq:secondary-residue} gives a meromorphic continuation of
the residue to a domain containing \(w=\frac{1}{2}\). In particular,
the value at the central point is obtained from this continued
expression.

We next identify this continuation with the residue of
\(A_3(s,w)\) on \(s+w=\frac{4}{3}\). To account for its
intersection with \(s=1\), restrict the holomorphic function
\eqref{eq:cleared} to \(s=\frac{4}{3}-w\). On the initial residue
domain, the functional equation \eqref{eq:FE} gives the value
\[
\left(\frac{1}{3}-w\right)\zeta_{c_K}(2)
 \pi^{w-\frac{1}{2}}
 \frac{\Gamma\left(\frac{1-w}{2}\right)}
      {\Gamma\left(\frac{w}{2}\right)}
 \Res_{u=\frac{4}{3}}C_{3,1}(u,1-w).
\]
In this expression, the pole at \(w=\frac{1}{3}\) comes from
\(\zeta\left(\frac{4}{3}-w\right)\) in
\eqref{eq:secondary-residue}; it is canceled by the factor
\(\frac{1}{3}-w\). The identity therefore continues through the
intersection and, by uniqueness, identifies the residue at
\(w=\frac{1}{2}\). At this point the gamma factor in
\eqref{eq:FE} equals \(1\), and the change of variable
\(s\mapsto s+w\) has derivative \(1\). Hence
\begin{equation}
\Res_{s=\frac{5}{6}}A_3\left(s,\frac{1}{2}\right)
=\Res_{s=\frac{4}{3}}C_{3,1}\left(s,\frac{1}{2}\right)
=C_3'.
\label{eq:central-residue}
\end{equation}

It remains to evaluate this residue and determine its sign.
Using \eqref{eq:kappap}, the local factor in the product in
\eqref{eq:secondary-residue} at \(w=\frac{1}{2}\) becomes
\[
\begin{cases}
\displaystyle
\frac{\left(1-p^{-\frac{5}{6}}\right)
\left(1+2p^{-1}+p^{-\frac{5}{6}}-2p^{-\frac{5}{2}}\right)}
{(1+p^{-1})^2},&p\equiv1\pmod3,\\[10pt]
\displaystyle
\frac{1-p^{-\frac{5}{3}}}{1+p^{-2}},&p\equiv2\pmod3.
\end{cases}
\]
This yields \eqref{eq:C3prime}. Each local factor is strictly
positive; for the split primes, note that
\(2p^{-1}-2p^{-\frac{5}{2}}>0\).
Moreover, both local factors are \(1+O(p^{-\frac{5}{3}})\), so their
deviations from \(1\) are summable over the primes. The two products
therefore converge to strictly positive values.
Finally, the alternating-series identity gives
\[
\zeta\left(\frac{5}{6}\right)
=\frac{\displaystyle\sum_{n=1}^{\infty}
(-1)^{n-1}n^{-\frac{5}{6}}}{1-2^{\frac{1}{6}}}<0,
\]
since the numerator is positive and the denominator is negative.
Every other factor in \eqref{eq:C3prime} is strictly positive.
Hence \(C_3'<0\), proving the asserted nonvanishing.

\begin{remark}
The normalization can also be checked against the first residue at
the intersection \((s,w)=\left(1,\frac{1}{3}\right)\). Multiplying the correction
product in \eqref{eq:secondary-residue} at \(w=\frac{2}{3}\) by \(\zeta(2)\)
gives exactly the product in \eqref{eq:first-residue} at \(w=\frac{1}{3}\).
At the ramified prime, the factor \(1-3^{-2}\) cancels the local
factor of \(\zeta(2)\). The gamma duplication formula gives
\[
\pi^{-\frac{1}{6}}\frac{\Gamma\left(\frac{1}{3}\right)}{\Gamma\left(\frac{1}{6}\right)}
 \frac{(2\pi)^{\frac{5}{3}}}
 {8\,3^{\frac{3}{2}}\Gamma\left(\frac{2}{3}\right)\zeta_K(2)}
=\frac{r_K}{4\zeta_K(2)}.
\]
Thus the two residue functions have opposite residues at \(w=\frac{1}{3}\),
as required by the denominator \((s-1)\left(s+w-\frac{4}{3}\right)\).
\end{remark}

\subsection{The secondary coefficient in the second moment}
\label{sec:second-moment}

We now prove Theorem~\ref{thm:second-moment}. The proof in this subsection is unconditional.
\allowdisplaybreaks[1]


\subsubsection{The natural family and normalizations}
\label{cub:setup}

We retain \(K,\omega,c_K\) from the introduction and put
\[
 \mathcal O_K=\Z[\omega],\qquad
 \lambda=c_K=1+2\omega=\sqrt{-3},\qquad
 S=S_\infty\cup\{(\lambda)\}.
\]
We fix this choice of \(\lambda\) throughout the second-moment calculation.
Although \((\lambda)=(1-\omega)\), the supplementary reciprocity laws
depend on the chosen generator. For any finite enlargement \(T\) of
\(S\), let \(I(T)\) denote the nonzero integral ideals coprime to the
finite primes in \(T\). The superscript \(T\) on an \(L\)-function
or on \(\zeta_K\) means that its Euler factors at those finite primes
have been removed. We use \(\mu_K\) for the ideal M\"obius function.

Each ideal in \(I(S)\) has a unique primary generator, and we identify
these ideals with their generators. Throughout the second-moment
calculation, primes are prime ideals of \(K\), represented by primary elements.
The second moment contains all squarefree ideals in \(I(S)\).
Its natural Hecke character is determined by
\[
 \chi_a((b))=\left(\frac{a}{b}\right)_3
 \qquad\bigl(b\text{ primary},\ (b,3a)=1\bigr).
\]
By cubic reciprocity this equals \(\chi_{3,a}(b)\), in the notation
of Section~\ref{sec:preliminaries}. On rational integers \(m\) with \((m,3a)=1\),
these values also agree with \(\widehat\chi_{3,a}(m)\), the Dirichlet
character used in the first moment. In functional equations, \(\chi_a\) always
denotes the associated primitive Hecke character on ideals of \(K\).
For \(\Re w\) sufficiently large, the squarefree series to be studied is
\begin{equation}\label{cub:squarefree-series-definition}
 Z_1^S\left(\frac12,\frac12,w\right)
 =\sum_{a\in I(S)}
   \frac{\mu_K(a)^2\left|L^S\left(\frac12,\chi_a\right)\right|^2}
        {N(a)^w}.
\end{equation}

For supplementary characters arising from units of
\(\mathcal O_{K,S}\), we use the separate notation
\begin{equation}\label{cub:supplementary-character}
 \xi_t(c)=\left(\frac{t}{c}\right)_3,
 \qquad t\in\mathcal O_{K,S}^{\times},\quad c\equiv1\pmod3.
\end{equation}
This is a finite character of the modulus variable \(c\).
It does not denote an arbitrary extension of the ideal character
\(\chi_c\) to the element \(t\). For example, the supplementary
value at \(\omega\) must be distinguished from the Hecke character
value at the unit ideal \((\omega)=(1)\).

Write \(e(x)=e^{2\pi i x}\). For
\(r\in\mathcal O_K\setminus\{0\}\), we use the Gauss sums
normalized by the trace in \cite[Section~3]{DDFDS}:
\begin{equation}\label{cub:raw-gauss-definition}
 g_3(r,c)=\sum_{x\bmod c}\left(\frac{x}{c}\right)_3
       e\!\left(\Tr\left(\frac{rx}{c}\right)\right),
 \qquad c\equiv1\pmod3,
 \qquad g_3(r,1)=1.
\end{equation}
The relation with \eqref{eq:gauss} is exactly
\[
 g_3(r,c)=g_K(c_Kr,\chi_{3,c}).
\]
If \((t,c)=1\), a change of variables gives
\begin{equation}\label{cub:raw-gauss-scaling}
 g_3(tr,c)=\overline{\left(\frac{t}{c}\right)_3}\,g_3(r,c).
\end{equation}
For coprime primary elements \(c,d\), the Chinese remainder theorem
and cubic reciprocity give
\begin{equation}\label{cub:raw-gauss-product}
 \begin{split}
 g_3(r,cd)
 &=\left(\frac{c}{d}\right)_3\left(\frac{d}{c}\right)_3
       g_3(r,c)g_3(r,d)\\
 &=\overline{\left(\frac{c}{d}\right)_3}\,
       g_3(r,c)g_3(r,d).
 \end{split}
\end{equation}
We distinguish these Gauss sums and their residues from the
normalized quantities in \cite{DIPP}. Put
\begin{equation}\label{cub:raw-series}
 \psi(r,v)=\sum_{c\equiv1\ (3)}
                \frac{g_3(r,c)}{N(c)^v},
 \qquad
 \rho(r)=\Res_{v=\frac43}\psi(r,v).
\end{equation}
Thus \(\psi\) is the unrestricted series already used in the proof
of Lemma~\ref{lem:Hbound}. We also write
\[
 \zeta_\lambda(s)=\zeta_{c_K}(s)=(1-3^{-s})\zeta_K(s),\qquad
 \kappa_\lambda=\Res_{s=1}\zeta_\lambda(s),\qquad
 u_0=\frac16.
\]
With the constant \(c_0\) already used in the proof of
Lemma~\ref{lem:Hresidue}, we put
\begin{equation}\label{cub:rho-constant}
 c_0=\frac{(2\pi)^{\frac53}}
 {8\cdot3^{\frac92}\Gamma\left(\frac23\right)\zeta_K(2)},
 \qquad \rho_0=27c_0>0.
\end{equation}

\subsubsection{Patterson residues and the removal of modulus primes}
\label{cub:raw-and-removal}

\begin{lemma}[Patterson residues]\label{cub:patterson-lemma}
Let \(a,b\) be primary, with \(a\) squarefree. Then
\begin{equation}\label{cub:raw-residue}
 \rho(ab^3)=\rho_0\,
       \frac{\overline{g_3(1,a)}}{N(a)^{\frac23}}.
\end{equation}
No coprimality between \(a\) and \(b\) is required. If the cubefree
part of a primary index contains a prime to exponent \(2\), its
residue is zero. In particular, for every primary prime \(p\) coprime
to \(3\) and every nonzero index \(r\),
\begin{equation}\label{cub:rho-periodicity}
 \rho(p^3r)=\rho(r).
\end{equation}
For primary \(a\), the nine shifts by units of \(\mathcal O_{K,S}\)
satisfy
\begin{equation}\label{cub:raw-unit-shifts}
 \rho(\omega^j a)=
 \begin{cases}
   \rho(a),&j=0,\\
   0,&j=1,2,
 \end{cases}
 \qquad
 \rho(\omega^j\lambda a)=0,
 \quad 0\leq j<3,
\end{equation}
and
\begin{equation}\label{cub:raw-lambda-shifts}
 \begin{split}
 \rho(\omega^j\lambda^2a)
 &=3^{-\frac13}e_j\,\xi_{\omega^j\lambda^2}(a)\rho(a),\\
 (e_0,e_1,e_2)&=\left(1,e\left(-\frac19\right),e\left(\frac19\right)\right).
 \end{split}
\end{equation}
These identities remain valid when the index is multiplied by a cube.
\end{lemma}

\begin{proof}
By \cite[Lemma~6.3]{DDFDS}, Patterson's residue is
\[
 \rho(r)=c_0\,\frac{\tau_3(r)}{N(r)^{\frac16}}.
\]
For \(r=ab^3\), the nonzero case of \cite[(6.8)]{DDFDS} gives
\[
 \tau_3(ab^3)
 =27\,\overline{g_3(1,a)}
       \frac{N(b)^{\frac12}}{N(a)^{\frac12}}.
\]
Since \(N(ab^3)^{\frac16}=N(a)^{\frac16}N(b)^{\frac12}\),
division by this factor proves \eqref{cub:raw-residue}.
The same formula for the theta coefficients vanishes when the cubefree
part contains a prime to exponent \(2\). Increasing any prime exponent
by \(3\) therefore leaves the residue unchanged, including when that
prime already divides the index. This proves
\eqref{cub:rho-periodicity}.

The remaining support conditions in \cite[(6.8)]{DDFDS} give
\eqref{cub:raw-unit-shifts}. For \(\omega^j\lambda^2ab^3\), the
nonzero theta coefficient is
\[
 27e_j\,\overline{g_3(\omega^j\lambda^2,a)}
             \frac{N(b)^{\frac12}}{N(a)^{\frac12}}.
\]
Relative to the index \(ab^3\), the sixth root of its norm acquires
the factor \(3^{\frac13}\). Formula
\eqref{cub:raw-gauss-scaling} now gives
\eqref{cub:raw-lambda-shifts}. Since \(\xi_t(b^3)=1\), its right
side may equally be written using \(ab^3\) in place of \(a\).
If the index contains a cube of \(\lambda\), the powers of \(3\)
in \cite[(6.8)]{DDFDS} and in \(N(r)^{\frac16}\) change by the same
factor, so the assertion also covers these cubes.
\end{proof}

The lemma and \eqref{cub:raw-gauss-product} imply, for \(p\nmid r\)
and \(q=N(p)\),
\begin{equation}\label{cub:raw-prime-recurrence}
 \rho(rp)=q^{-\frac23}\overline{g_3(r,p)}\,\rho(r),
 \qquad \rho(rp^2)=0.
\end{equation}
Here \(r\) may contain factors of \(\omega\) or \(\lambda\).
In the nonzero case of \eqref{cub:raw-lambda-shifts}, adjoining \(p\)
contributes the supplementary factor
\(\left(\frac{\omega^j\lambda^2}{p}\right)_3\), which is already
included in \(\overline{g_3(r,p)}\). In the other cases both sides
vanish.

\begin{lemma}[Removing primes from the moduli]\label{cub:modulus-removal-lemma}
For a squarefree primary \(h\) coprime to \(3\), define
\[
 \psi_h(r,v)=
 \sum_{\substack{c\equiv1\ (3)\\(c,h)=1}}
       \frac{g_3(r,c)}{N(c)^v},
 \qquad
 \rho_h(r)=\Res_{v=\frac43}\psi_h(r,v).
\]
Then, for every nonzero index \(r\),
\begin{equation}\label{cub:modulus-removal}
 \rho_h(r)=\beta_h\rho(r),
 \qquad
 \beta_h=\prod_{p\mid h}\left(1+N(p)^{-1}\right)^{-1}.
\end{equation}
In particular, \eqref{cub:modulus-removal} also holds when
\((r,h)\ne1\).
\end{lemma}

\begin{proof}
We induct on the number of prime divisors of \(h\). Write
\(h=\gamma p\), with \((p,\gamma)=1\), and
\(r=\delta p^k\), where \(p\nmid\delta\). Put
\(q=N(p)\) and \(\Delta_p(v)=1-q^{2-3v}\).
Expanding according to the exponent of \(p\) in the modulus gives
the local identities \cite[(6.5)--(6.7)]{DDFDS}:
\begin{align}
 \Delta_p(v)\psi_{\gamma p}(\delta,v)
 &=\psi_\gamma(\delta,v)
      -q^{-v}g_3(\delta,p)\psi_\gamma(\delta p,v),
       \label{cub:remove-zero}\\
 \Delta_p(v)\psi_{\gamma p}(\delta p,v)
 &=\psi_\gamma(\delta p,v)
      -q^{1-2v}\overline{g_3(\delta,p)}\psi_\gamma(\delta,v),
       \label{cub:remove-one}\\
 \Delta_p(v)\psi_{\gamma p}(\delta p^2,v)
 &=\psi_\gamma(\delta p^2,v).
       \label{cub:remove-two}
\end{align}
Take residues at \(v=\frac43\). The induction hypothesis,
\eqref{cub:raw-prime-recurrence}, and
\(|g_3(\delta,p)|^2=q\) turn the right side of
\eqref{cub:remove-zero} into
\[
 \beta_\gamma\rho(\delta)
 \left(1-q^{-\frac43}q^{-\frac23}|g_3(\delta,p)|^2\right)
 =\beta_\gamma(1-q^{-1})\rho(\delta).
\]
Similarly, the right side of \eqref{cub:remove-one} becomes
\(\beta_\gamma(1-q^{-1})\rho(\delta p)\).
Since \(\Delta_p(\frac43)=1-q^{-2}\), both exponent classes
give the factor
\[
 \frac{1-q^{-1}}{1-q^{-2}}=\frac1{1+q^{-1}}.
\]
The class \(k=2\) follows from \eqref{cub:remove-two} and
\(\rho(\delta p^2)=0\).

Every modulus in \(\psi_{\gamma p}\) is coprime to \(p\), and
\eqref{cub:raw-gauss-scaling} gives
\[
 g_3(p^3r,c)=g_3(r,c).
\]
Thus \(\psi_{\gamma p}(p^3r,v)=\psi_{\gamma p}(r,v)\), first
term by term in the region of absolute convergence and then by
continuation. Together with \eqref{cub:rho-periodicity}, this reduces
every \(k\geq0\) to the three cases above and completes the induction.
\end{proof}

\begin{proposition}[The square sum of Patterson residues]
\label{cub:theta-square-proposition}
For \(\Re u>\frac23\), let
\[
 T_3(u)=\sum_{a\equiv1\ (3)}\frac{|\rho(a)|^2}{N(a)^u}.
\]
Then
\begin{equation}\label{cub:theta-square-sum}
 \begin{split}
 T_3(u)
 &=\rho_0^2\prod_{p\ne\lambda}
       \frac{1+q^{-u-\frac13}}{1-q^{-3u}}\\
 &=\rho_0^2
   \frac{\zeta_\lambda(3u)\zeta_\lambda\left(u+\frac13\right)}
        {\zeta_\lambda\left(2u+\frac23\right)},
 \qquad q=N(p).
 \end{split}
\end{equation}
This quotient gives its meromorphic continuation, and
\begin{equation}\label{cub:theta-derivative}
 T_3(u_0)=0,\qquad
 T_3'(u_0)=
   \frac{2\rho_0^2\zeta_\lambda\left(\frac12\right)^2}
        {\kappa_\lambda}>0.
\end{equation}
\end{proposition}

\begin{proof}
Every index with nonzero residue has a unique expression \(ab^3\),
with \(a\) squarefree; \(a\) and \(b\) may have common factors.
By \eqref{cub:raw-residue} and \(|g_3(1,a)|^2=N(a)\),
\[
 |\rho(ab^3)|^2=\rho_0^2N(a)^{-\frac13}.
\]
At a prime of norm \(q\), the exponents \(3m,3m+1,3m+2\)
therefore contribute \(q^{-3mu}\), \(q^{-3mu-u-\frac13}\), and
\(0\), respectively. Summing them gives
\eqref{cub:theta-square-sum}. Near \(u=u_0\),
\[
 \zeta_\lambda\left(2u+\frac23\right)
 =\frac{\kappa_\lambda}{2(u-u_0)}+O(1),
\]
whereas the numerator in the quotient has value
\(\zeta_\lambda(\frac12)^2\).

To verify that this numerator is nonzero, use
\(\zeta_K(s)=\zeta(s)L(s,\chi_{-3})\). For real \(0<s<1\),
the alternating series gives
\[
 \zeta(s)=
 \frac{\sum_{n\geq1}(-1)^{n-1}n^{-s}}{1-2^{1-s}}<0,
\]
while
\[
 L(s,\chi_{-3})
 =\sum_{n\geq0}
   \left((3n+1)^{-s}-(3n+2)^{-s}\right)>0.
\]
The grouped series in the last display converges and every term is
positive. Hence \(\zeta_\lambda(\frac12)\ne0\), and the expansion
above proves \eqref{cub:theta-derivative} without a hypothesis on zeros.
\end{proof}

\subsubsection{Gauss sums and finite characters}
\label{cub:gauss-normalization}

In the construction of \cite{DIPP}, take the ray modulus
\(\mathfrak c=(9)=(\lambda^4)\). We first identify its full finite
class group. Put
\[
 U_n=1+\lambda^n\mathcal O_{K,\lambda}.
\]
The local unit group is \(\mathcal O_{K,\lambda}^{\times}=\mu_6U_2\),
with \(\mu_6\cap U_2=\{1\}\). The local logarithm identifies
\(U_2\) with \(\lambda^2\mathcal O_{K,\lambda}\). Since
\(3=-\lambda^2\),
\[
 \log(U_2^3)=3\lambda^2\mathcal O_{K,\lambda}
           =\lambda^4\mathcal O_{K,\lambda}
           =\log(U_4),
 \qquad U_2^3=U_4.
\]
This also verifies the local cube condition on \(\mathfrak c\) in
\cite{DIPP}. Since \(K\) has class number one and global unit group
\(\mu_6\), its ray class group modulo \(\mathfrak c\) is
\[
 (\mathcal O_K/(\lambda^4))^\times/\mu_6
 \cong U_2/U_4.
\]
This group has order \(9\) and exponent \(3\). Every finite Hecke
character of order dividing \(3\), unramified outside \(\lambda\),
is trivial on \(\mu_6\) and on local cubes, and therefore factors
through \(U_2/U_4\). Conversely, every character of this quotient
defines such a ray class character. Thus the full group
\(R_{\mathfrak c}\) of \cite[Section~2]{DIPP} and its dual are
precisely this group and its dual. Write
\[
 \mathcal E=R_{\mathfrak c}
 \cong U_2/U_4\cong\mathbb F_3^2,
\]
and represent its elements by primary residue classes modulo \(9\).
If
\[
 c\equiv1+\alpha_2\lambda^2+\alpha_3\lambda^3\pmod9,
 \qquad \alpha_2,\alpha_3\in\{-1,0,1\},
\]
the supplementary laws \cite[(3.2)]{DDFDS} give
\begin{equation}\label{cub:class-coordinates}
 \left(\frac{\omega}{c}\right)_3=\omega^{\alpha_2},
 \qquad
 \left(\frac{\lambda}{c}\right)_3=\omega^{-\alpha_3}.
\end{equation}
We use coordinates \(E=(j,k)\) determined by
\[
 \xi_\omega(E)=\omega^j,\qquad \xi_\lambda(E)=\omega^k.
\]
The coordinates are independent, and multiplication in \(\mathcal E\)
adds them. In particular, the cube of a primary element represents
the identity class.

Choose primary generators for the class representatives in
\cite{DIPP}. Cubic reciprocity then makes the additional reciprocity
factor \(\psi_E(a)\) equal to \(1\) on primary indices \(a\).
Indeed, global reciprocity identifies the defining product of local
Hilbert symbols at \(S\) with
\[
 \frac{\left(\frac{a}{m_E}\right)_3}
      {\left(\frac{m_E}{a}\right)_3}=1.
\]
For noncoprime parameters one first chooses coprime representatives
in the same classes. Thus these reciprocity factors introduce no
additional character of the primary outer index.

\begin{lemma}[The finite phase of the Gauss sum]\label{cub:gauss-phase-lemma}
Define
\begin{equation}\label{cub:gauss-phase}
 f_\lambda(j,k)=
 \begin{cases}
   \omega^k,&j=0,\\
   1,&j=1,2.
 \end{cases}
\end{equation}
The generalized Gauss sums normalized as in
\cite[(5), Section~4.1.1]{DIPP} satisfy
\begin{equation}\label{cub:generalized-gauss-match}
 G(a,b)=f_\lambda([b])\,\frac{g_3(a,b)}{N(b)^{\frac12}}
\end{equation}
for all primary \(a,b\), including nonsquarefree moduli \(b\).

For the associated primitive root number, let \(b_0\) be the cubefree
part of \(b\), and put
\[
 c=\operatorname{rad}(b_0),\qquad
 \eta_b(x)=\left(\frac{x}{b_0}\right)_3.
\]
Then
\begin{equation}\label{cub:primitive-gauss-match}
 G(\chi_b)=
 f_\lambda([b])\,N(c)^{-\frac12}
 \sum_{x\bmod c}\eta_b(x)
    e\!\left(\Tr\left(\frac{x}{c}\right)\right).
\end{equation}
When a prime occurs to exponent \(2\) in \(b_0\), this expression
uses a Gauss sum modulo \(c\), and cannot be replaced by
\(f_\lambda([b])\frac{g_3(1,b_0)}{N(c)^{\frac12}}\).
\end{lemma}

\begin{proof}
We begin with the primitive local factor at \(\lambda\). Write
\([b]=(j,k)\). Since \(\eta_b(\omega)=\omega^j\), the local
factor needed to make the residue class character trivial on global
units is \(\theta^j\), where
\[
 \theta(\omega)=\omega^{-1},\qquad
 \theta(-1)=1,\qquad
 \theta|_{1+\lambda^2\mathcal O_{K,\lambda}}=1.
\]
Its conductor exponent is
\[
 \nu_j=
 \begin{cases}
   0,&j=0,\\
   2,&j=1,2.
 \end{cases}
\]
Thus the primitive residue class character defining \(\chi_b\)
is \(\eta_b\theta^j\), of conductor \(\lambda^{\nu_j}c\).
For \(j=1,2\), its finite local Gauss factor is
\[
 \tau_j=\frac13\sum_{x\bmod\lambda^2}^{\times}
       \theta(x)^j
       e\!\left(\Tr\left(\frac{x}{\lambda^3}\right)\right).
\]
Use the representatives \(\{\pm1,\pm\omega,\pm\omega^2\}\).
The sums of the additive character on these three pairs are
\(2,-1,-1\), while the values of \(\theta^j\) are
\(1,\omega^{-j},\omega^{-2j}\). Consequently,
\begin{equation}\label{cub:local-gauss-sum}
 \tau_j=\frac{2-\omega^{-j}-\omega^{-2j}}3=1,
 \qquad j=1,2.
\end{equation}
For \(j=0\), the conductor is \(1\), and its local Gauss factor
is \(\tau_0=1\). The sum of conductor \(\lambda^2\) above is not
used in this case.

The additive character in \cite{DIPP} uses the different generator
\(\lambda\). In its primitive Gauss sum, write
\(x=cy+\lambda^{\nu_j}z\) by the Chinese remainder theorem and then
substitute \(z=\lambda z'\) in the sum modulo \(c\). This gives
\[
 G(\chi_b)
 =\tau_j\,\eta_b(\lambda^{\nu_j+1})\,N(c)^{-\frac12}
   \sum_{z'\bmod c}\eta_b(z')
        e\!\left(\Tr\left(\frac{z'}c\right)\right).
\]
Since \(c\) is primary, \(\theta(c)=1\), and
\(\eta_b(\lambda)=\omega^k\). The finite phase is therefore
\(\omega^k\) when \(j=0\), and \(\eta_b(\lambda^3)=1\)
when \(j\ne0\). This proves \eqref{cub:primitive-gauss-match}.

We next check every prime power in the modulus. Fix \(p\ne\lambda\)
and put
\[
 q=N(p),\qquad k=v_p(a),\qquad \ell=v_p(b),\qquad
 a_p=\frac{a}{p^k},\qquad
 \eta_p(x)=\left(\frac{x}{p}\right)_3.
\]
Write
\[
 \mathfrak g(\eta_p^\ell;p)
 =\sum_{x\bmod p}\eta_p(x)^\ell
       e\!\left(\Tr\left(\frac{x}{p}\right)\right).
\]
Summing first over residue classes modulo \(p\) yields
\begin{equation}\label{cub:prime-power-table}
 \frac{g_3(a,p^\ell)}{q^{\frac\ell2}}=
 \begin{cases}
  1,&\ell=0,\\
  q^{\frac{\ell-2}{2}}\overline{\eta_p(a_p)^\ell}
        \mathfrak g(\eta_p^\ell;p),
      &k=\ell-1,\quad 3\nmid\ell,\\
  -q^{\frac{\ell-2}{2}},
      &k=\ell-1,\quad \ell>0,\quad3\mid\ell,\\
  q^{\frac\ell2}(1-q^{-1}),
      &k\geq\ell,\quad \ell>0,\quad3\mid\ell,\\
  0,&\text{otherwise}.
 \end{cases}
\end{equation}
In the second line, the sum modulo \(p\) is repeated
\(q^{\ell-1}\) times. The third and fourth lines come from the
Ramanujan sums \(-q^{\ell-1}\) and \(q^\ell-q^{\ell-1}\),
respectively, followed by division by \(q^{\frac\ell2}\).
These are the local factors \(G_0(p^k,p^\ell)\) in
\cite[Section~4.1.1]{DIPP}, with the appropriate normalized primitive
Gauss sum.

It remains to match the factors involving distinct primes. Let
\(P=\{p:3\nmid v_p(b)\}\). For a nonzero local term,
\(k_p=\ell_p-1\) at every \(p\in P\). For \(p,q\in P\), combine
the Chinese remainder factor for \(g_3(a,b)\) with the unit part of
the index in \eqref{cub:prime-power-table}. The relations
\[
 \ell_q-k_q=1,\qquad \ell_p-k_p=1
\]
reduce their product to
\(\eta_p(q)^{\ell_p}\eta_q(p)^{\ell_q}\). This is the Chinese
remainder factor in the primitive Gauss sum of modulus
\(\operatorname{rad}(b_0)\). The remaining index factors combine to
\(\overline{\eta_b(a^{(b_0)})}\), where \(a^{(b_0)}\) is the part
of \(a\) coprime to \(b_0\). If the exponent of a modulus prime is
divisible by \(3\), its modulus character is trivial, and its index
contribution is already included in
\(\overline{\eta_b(a^{(b_0)})}\). Hence
\eqref{cub:prime-power-table} and \eqref{cub:primitive-gauss-match}
give precisely the definition
\[
 G(a,b)=
 \overline{\chi_{b_0}(a^{(b_0)})}\,
 G(\chi_{b_0})G_0(a,b)
\]
from \cite{DIPP}. Formula \eqref{cub:generalized-gauss-match} follows.
\end{proof}

\begin{lemma}[The quotient of the finite phases]
\label{cub:gamma-cocycle-lemma}
Let \(\gamma(E)=\gamma(\psi_E)\), with the definitions in
\cite[(6)--(7), Lemma~3.3.2]{DIPP}. Then
\begin{equation}\label{cub:gamma-cocycle}
 \vartheta(E)=\frac{\gamma(E)}{f_\lambda(E)}
\end{equation}
is a character of \(\mathcal E\). In particular, for
\(E,E',F\in\mathcal E\),
\begin{equation}\label{cub:gamma-cancellation}
 \begin{split}
 &\gamma(E)\overline{\gamma(E')}\gamma(F)
       \gamma(E'E^{-1}F^{-1})\\
 &\qquad =
 f_\lambda(E)\overline{f_\lambda(E')}f_\lambda(F)
       f_\lambda(E'E^{-1}F^{-1}).
 \end{split}
\end{equation}
\end{lemma}

\begin{proof}
Choose coprime cubefree primary parameters representing \(E_1,E_2\).
The root number multiplication law in \cite[Lemma~3.3.2]{DIPP} has
the cocycle
\[
 \tau(E_1,E_2)=
 \frac{\gamma(E_1E_2)}{\gamma(E_1)\gamma(E_2)}.
\]
On the other hand, apply the Chinese remainder theorem to the two
coprime primitive moduli in \eqref{cub:primitive-gauss-match}.
The contribution from the primes outside \(S\) is
\[
 \chi_{b_1}(\operatorname{rad}b_2)
 \chi_{b_2}(\operatorname{rad}b_1).
\]
These are exactly the two character values divided out in the
definition of \(\tau\) in \cite{DIPP}. The remaining cocycle is
therefore
\[
 \tau(E_1,E_2)=
 \frac{f_\lambda(E_1E_2)}
      {f_\lambda(E_1)f_\lambda(E_2)}.
\]
Dividing the two expressions shows that
\(\vartheta(E_1E_2)=\vartheta(E_1)\vartheta(E_2)\).
Every pair of classes admits coprime primary representatives, so
this identity holds on all of \(\mathcal E\). Since each phase has
absolute value one,
\[
 \vartheta(E)\overline{\vartheta(E')}\vartheta(F)
 \vartheta(E'E^{-1}F^{-1})=1,
\]
which proves \eqref{cub:gamma-cancellation}.
\end{proof}

\begin{remark}\label{cub:gamma-extension-remark}
Formula \eqref{cub:gamma-cancellation} concerns the full product of
four phases. It does not require a choice for which
\(\gamma(E)=f_\lambda(E)\) separately at every class. Thus the
calculation is independent of how the auxiliary characters
\(\psi_E\) are extended to units of \(\mathcal O_{K,S}\).
The Euler factor at \(\lambda\) below is nevertheless determined
by the actual primitive Hecke character.
\end{remark}

To compare the two residue normalizations, define
\begin{equation}\label{cub:class-residue-definition}
 v_E(a)=\Res_{v=\frac43}
 \sum_{\substack{c\equiv1\ (3)\\{[c]=E}}}
             \frac{g_3(a,c)}{N(c)^v}.
\end{equation}
Projection onto the nine finite characters expresses this residue
in terms of the preceding shifted series, and supplies its
meromorphic continuation. Let \(h\) be squarefree,
\(T=S\cup S_h\), and let \(f\) be supported on the prime divisors
of \(h\), where \(S_h\) is the set of those primes. For a finite
character \(\xi\) of \(\mathcal E\), the relevant class component
of the Kubota series is
\begin{equation}\label{cub:class-kubota-definition}
 D^T(s,a,\xi\chi_f;E)
 =\sum_{\substack{b\in I(T)\\{[b]=E}}}
    \frac{\xi(b)\chi_f(b)G(a,b)}{N(b)^s}.
\end{equation}
Its normalization in \cite{DIPP} is
\[
 \widetilde D^T(s,a,\xi\chi_f;E)
 =\zeta_K^T\left(3s-\frac12\right)D^T(s,a,\xi\chi_f;E).
\]
By \eqref{cub:generalized-gauss-match}, the identity
\(\chi_f(c)g_3(a,c)=g_3(f^2a,c)\), and
\eqref{cub:modulus-removal}, its residue is
\begin{equation}\label{cub:normalized-residue}
 \begin{split}
 \widetilde R_E^T(a;\xi\chi_f)
 &:=\Res_{s=\frac56}
       \widetilde D^T(s,a,\xi\chi_f;E)\\
 &=\zeta_K^T(2)\beta_h f_\lambda(E)\xi(E)v_E(f^2a)\\
 &=\zeta_\lambda(2)
     \prod_{p\mid h}(1-q^{-1})
       f_\lambda(E)\xi(E)v_E(f^2a),
 \qquad q=N(p).
 \end{split}
\end{equation}
The formula for removing modulus primes in each class follows by
finite Fourier projection from \eqref{cub:modulus-removal}, which
holds for every shifted index. It therefore remains valid when
\(f^2a\) and \(h\) have a common factor. The last equality uses
\[
 \frac{\zeta_K^T(2)}{\zeta_\lambda(2)}
      \prod_{p\mid h}(1+q^{-1})^{-1}
 =\prod_{p\mid h}(1-q^{-1}).
\]
The variables are related by \(v=s+\frac12\), so taking residues
introduces no further scalar factor.

\subsubsection{Root numbers and the Euler factor at \texorpdfstring{\(\lambda\)}{lambda}}
\label{cub:corrected-functional-equation}

\begin{proposition}[The root number for nontrivial \(g\)]
\label{cub:root-correction-proposition}
Use the natural primary normalization above, take the additional
finite character in the third functional equation of \cite{DIPP}
to be trivial, and retain \(g\mid h^2\). If
\(C_{\mathrm{printed}}(E,E',g,1,w)\) denotes the coefficient
displayed in \cite[Proposition~5.3.2]{DIPP}, the coefficient needed
here is
\begin{equation}\label{cub:corrected-root-number}
 C_{\mathrm{corr}}(E,E',g;w)
 =G(\overline{\chi_g})\,
     C_{\mathrm{printed}}(E,E',g,1,w).
\end{equation}
More explicitly, write
\(\Delta=E'E^{-1}[g]^{-1}=(j_\Delta,k_\Delta)\). Then
\begin{equation}\label{cub:corrected-class-coefficient}
 \begin{split}
 C_{\mathrm{corr}}(E,E',g;w)
 ={}&G(\overline{\chi_g})\,
 f_\lambda(E)\overline{f_\lambda(E')}
 f_\lambda([g])f_\lambda(\Delta)\\
 &{}\times
 3^{(1+\nu_{j_\Delta})(w-\frac12)},
 \qquad \nu_0=0,\quad \nu_1=\nu_2=2.
 \end{split}
\end{equation}
Here \(\nu_j\) is the local conductor exponent computed above.
\end{proposition}

\begin{proof}
The required root number can first be checked in the primitive Hecke
functional equation. Put
\[
 \Gamma_K(s)=(2\pi)^{-s}\Gamma(s),\qquad
 \mathfrak f_g=\operatorname{Cond}(\chi_g).
\]
Writing \(T_f\) for the finite primes in \(T\), we have
\begin{equation}\label{cub:primitive-functional-equation}
 \begin{split}
 L^T(1-w,\overline{\chi_g})
 ={}&G(\overline{\chi_g})\,
       (3N\mathfrak f_g)^{w-\frac12}
       \frac{\Gamma_K(w)}{\Gamma_K(1-w)}
       L^T(w,\chi_g)\\
 &{}\times\prod_{v\in T_f}
       \frac{L_v(w,\chi_g)}
            {L_v(1-w,\overline{\chi_g})}.
 \end{split}
\end{equation}
The term \(b=b'=1\) in the auxiliary series leaves precisely
\(L^T(1-w,\overline{\chi_g})\) as its outer sum. Transforming this
to \(L^T(w,\chi_g)\) therefore requires the factor
\(G(\overline{\chi_g})\). For general \(b,b'\), the same factor
remains throughout the Gauss sum multiplication in the proof of
\cite[Proposition~5.3.2, p.~24]{DIPP}. It appears explicitly in the
multiplication identity on that page; retaining it to the final
formula gives \eqref{cub:corrected-root-number}. Thus this is the
coefficient of the same functional equation, as verified directly
by the primitive functional equation.

With primary representatives, the factors \(\psi_{E'}(-1)\) and
\(\psi_{E'E^{-1}}(m_{[g_1]})\) in the displayed coefficient are
both \(1\). Replace the product of its four \(\gamma\) factors
using \eqref{cub:gamma-cancellation}. The conductor exponent at
\(\lambda\) is \(\nu_{j_\Delta}\), and \(|D_K|=3\).
These substitutions prove \eqref{cub:corrected-class-coefficient}.
\end{proof}

\begin{remark}\label{cub:root-correction-size}
The root number \(G(\overline{\chi_g})\) has absolute value one.
It therefore leaves the continuation and growth estimates in
\cite{DIPP} unchanged. Its phase contributes to the exact residue
and must be retained in the calculation of nonvanishing.
\end{remark}

For later use, write \(g=G_1G_2^2\), where \(G_1,G_2\) are
coprime and squarefree. Formula \eqref{cub:primitive-gauss-match}
then becomes
\begin{equation}\label{cub:general-g-root}
 \begin{split}
 G(\chi_g)
 &=f_\lambda([g])\,
   \frac{g_3(1,G_1)\overline{g_3(1,G_2)}}
        {N(G_1G_2)^{\frac12}},\\
 G(\overline{\chi_g})
 &=\overline{f_\lambda([g])}\,
   \frac{\overline{g_3(1,G_1)}g_3(1,G_2)}
        {N(G_1G_2)^{\frac12}}.
 \end{split}
\end{equation}
The Chinese remainder factor for the two conjugate characters is
\[
 \left(\frac{G_2}{G_1}\right)_3
 \left(\frac{G_1}{G_2}\right)_3^2=1.
\]
Also, \(\left(\frac{-1}{G_2}\right)_3=1\), so conjugating the Gauss
sum adds no sign. These are the full primitive root numbers that
will cancel the phases involving distinct sieve primes.

We now take \(h=f=g=1\) and compute the local contribution at
\(\lambda\) to the finite transformation. Let the two modulus
classes be \(E,E'\), and write
\[
 E'E^{-1}=(j,k),\qquad w=1-u.
\]
The normalized residues contribute \(f_\lambda(E)\) and
\(\overline{f_\lambda(E')}\). Division by
\eqref{cub:corrected-class-coefficient} cancels these two class
phases, leaving \(f_\lambda(j,k)^{-1}\) and the conductor factor
\(3^{(1+\nu_j)(u-\frac12)}\).

The third functional equation also requires the ratio of the local
Euler factors. In its expansion, let
\(\chi=\chi_m\overline{\chi_{m'}}\). Since
\([m'e][me]^{-1}=E'E^{-1}\), the parameter class of \(\chi\) is
\((-j,-k)\). When \(j=0\), this character is unramified at
\(\lambda\), with the actual primitive local values
\[
 \chi((\lambda))=\omega^{-k},\qquad
 \overline\chi((\lambda))=\omega^k.
\]
These values come from the primitive extension of
\(\eta_m\overline{\eta_{m'}}\). The removal of Euler factors in
\cite[Sections~5.3--5.4]{DIPP}, or directly
\eqref{cub:primitive-functional-equation}, gives the required ratio
\[
 \frac{L_\lambda(u,\overline\chi)}
      {L_\lambda(1-u,\chi)}
 =\frac{1-3^{u-1}\omega^{-k}}{1-3^{-u}\omega^k}.
\]
For \(j\ne0\), the conductor at \(\lambda\) is \(\lambda^2\)
and both Euler factors are \(1\). The full finite kernel is therefore
\begin{equation}\label{cub:lambda-kernel}
 H(j,k;u)=
 \begin{cases}
  \displaystyle
  3^{u-\frac12}\omega^{-k}
       \frac{1-3^{u-1}\omega^{-k}}{1-3^{-u}\omega^k},
       &j=0,\\[6pt]
  27^{u-\frac12},&j=1,2.
 \end{cases}
\end{equation}
In particular, \(j=0\) specifies the unramified case, but does not
imply \(\chi((\lambda))=1\). The nontrivial unramified values are
included in \eqref{cub:lambda-kernel}.

\subsubsection{The finite Fourier pairing}
\label{cub:finite-pairing}

For \(v,t\in\mathbb F_3\), define
\begin{equation}\label{cub:finite-fourier-definition}
 \xi_{v,t}(E)=\omega^{vj+tk}\quad(E=(j,k)),
 \qquad
 \mathcal R_{v,t}(a)=\sum_{E\in\mathcal E}\xi_{v,t}(E)v_E(a).
\end{equation}
Formula \eqref{cub:raw-gauss-scaling} gives
\begin{equation}\label{cub:fourier-as-shift}
 \mathcal R_{v,t}(a)=\rho(\omega^{2v}\lambda^{2t}a).
\end{equation}
Consequently, \eqref{cub:raw-unit-shifts} and
\eqref{cub:raw-lambda-shifts} determine the entire Fourier support:
\begin{equation}\label{cub:fourier-support}
 \begin{split}
 |\mathcal R_{0,0}(a)|^2&=|\rho(a)|^2,\\
 |\mathcal R_{v,1}(a)|^2&=3^{-\frac23}|\rho(a)|^2
       \quad(v=0,1,2).
 \end{split}
\end{equation}
The other five components vanish.

\begin{proposition}[Evaluation of the finite pairing]
\label{cub:finite-pairing-proposition}
For every primary index \(a\), the kernel
\eqref{cub:lambda-kernel} satisfies
\begin{equation}\label{cub:lambda-pairing}
 \sum_{E,E'\in\mathcal E}
 H(E'E^{-1};u)v_E(a)\overline{v_{E'}(a)}
 =d_\lambda(u)|\rho(a)|^2,
\end{equation}
where
\begin{equation}\label{cub:lambda-scalar}
 d_\lambda(u)=
 \frac{2\cdot3^{u-\frac12}}{9(1-3^{-3u})}
       \left(3^{-u}+3^{\frac13}\right)
       +\frac23\,27^{u-\frac12}.
\end{equation}
This function is holomorphic near \(u_0=\frac16\), and
\begin{equation}\label{cub:lambda-value}
 d_\lambda(u_0)=
 \frac{4}{9(1-3^{-\frac12})}
 =\frac23(1+3^{-\frac12})>0.
\end{equation}
\end{proposition}

\begin{proof}
Put
\[
 y=3^{-u},\quad z=3^{u-1},\quad
 A=3^{u-\frac12},\quad B=27^{u-\frac12}.
\]
Since \(yz=\frac13\), the identity
\[
 \frac1{1-y\omega^k}
 =\frac{1+y\omega^k+y^2\omega^{2k}}{1-y^3}
\]
gives
\begin{equation}\label{cub:kernel-fourier-expansion}
 H(0,k;u)=
 \frac{A}{1-y^3}
 \left\{\frac23y+(y^2-z)\omega^k
                   +\frac23\omega^{2k}\right\}.
\end{equation}
Denote the three coefficients by
\[
 a_0=\frac{2Ay}{3(1-y^3)},\qquad
 a_1=\frac{A(y^2-z)}{1-y^3},\qquad
 a_2=\frac{2A}{3(1-y^3)}.
\]
With the Fourier convention in
\eqref{cub:finite-fourier-definition}, set
\[
 \widehat H(v,t;u)
 =\sum_{j,k\in\mathbb F_3}
      H(j,k;u)\omega^{vj+tk}.
\]
Orthogonality in the \(k\) variable yields
\begin{equation}\label{cub:kernel-fourier-values}
 \widehat H(0,0;u)=3a_0+6B,\qquad
 \widehat H(v,1;u)=3a_2\quad(v=0,1,2).
\end{equation}
For the second identity, the two rows \(j=1,2\) contribute zero
because they are constant in \(k\). The remaining Fourier values
multiply zero residue components by \eqref{cub:fourier-support}.

Finite Fourier inversion gives
\[
 v_E(a)=\frac19
       \sum_{v,t\in\mathbb F_3}
          \overline{\xi_{v,t}(E)}\mathcal R_{v,t}(a).
\]
Substitute this expression in the double sum and first sum over
\(E\). The result is
\begin{equation}\label{cub:fourier-pairing-identity}
 \sum_{E,E'}H(E'E^{-1};u)v_E(a)\overline{v_{E'}(a)}
 =\frac19\sum_{v,t}
       \widehat H(v,t;u)|\mathcal R_{v,t}(a)|^2.
\end{equation}
The factor \(\frac19\) consists of \(\frac1{81}\) from the two
inversions and a factor \(9\) from character orthogonality.
Using \eqref{cub:fourier-support} and
\eqref{cub:kernel-fourier-values}, the coefficient of
\(|\rho(a)|^2\) on the right is
\[
 \frac19\left(3a_0+6B
             +3\cdot3a_2\cdot3^{-\frac23}\right)
 =\frac{a_0}{3}+\frac{2B}{3}+3^{-\frac23}a_2.
\]
Substitution proves \eqref{cub:lambda-scalar}. The calculation also
holds when \(\rho(a)=0\), since no division by this value is needed.
Finally, put \(r=3^{-\frac12}\). At \(u=u_0\),
\[
 d_\lambda(u_0)=\frac{2(1+r)}{9(1-r)}+\frac29
              =\frac4{9(1-r)}
              =\frac23(1+r),
\]
which proves \eqref{cub:lambda-value}.
\end{proof}

\subsubsection{The complete residue and its translations}
\label{cub:complete-pairing}

Let \(\widetilde Z^S\) denote the perfect multiple Dirichlet series
of \cite[Sections~5.5--5.6]{DIPP}, defined explicitly in
\eqref{cub:sieve-perfect-definition} below. Write \(\mathcal R_3(u)\) for
the complete pairing obtained by taking the two Kubota residues in
its third functional equation, including the finite class projections,
Gauss phases, and local Euler factors. Its normalization is fixed by
\begin{equation}\label{cub:complete-residue-normalization}
 \begin{split}
 &\lim_{s\to u+\frac13}
    \left(s-u-\frac13\right)^2
    \widetilde Z^S(s,s,1-u;1,1,1)\\
 &\qquad=
 \zeta_\lambda\left(2u+\frac23\right)
 \frac{\Gamma_K(u)}{\Gamma_K(1-u)}\mathcal R_3(u).
 \end{split}
\end{equation}
This is \cite[Proposition~5.7.2]{DIPP} with the parameters
\(h=f=g=1\). The corresponding identity for each sieve term is
\eqref{cub:sieve-double-residue-identity} below.

In the class basis, \eqref{cub:normalized-residue} supplies two
factors of \(\zeta_\lambda(2)\). Formulas
\eqref{cub:gamma-cancellation} and
\eqref{cub:corrected-class-coefficient} cancel the two class phases,
and the remaining transformation is the kernel
\eqref{cub:lambda-kernel}.

To justify taking the two residues term by term, first work in a
sufficiently far right half plane. Fix \(0<\delta<\frac1{12}\).
On the circle \(\lvert v-\frac43\rvert=\delta\),
\cite[Lemma~6.3]{DDFDS} gives the uniform bound
\[
 |\psi(r,v)|
 \ll_{\delta,\eps}
 N(r)^{\frac1{12}+\frac\delta2+\eps}.
\]
Indeed, this is the largest value of the exponent
\(\frac34-\frac{\Re v}{2}+\eps\) on that circle. The finite set
of shifts by units of \(\mathcal O_{K,S}\) only changes the implied
constant. On the corresponding circle
\(\lvert s-\frac56\rvert=\delta\), the factor
\(\zeta_\lambda(3s-\frac12)\) is also uniformly bounded.
Consequently, for \(\Re u>\frac76+\delta+2\eps\), the outer
\(a\)-series of the product of two Kubota factors converges
absolutely and uniformly on the product of the two circles.
Cauchy's formula therefore permits the two residue integrals to be
interchanged with the outer sum. For fixed \(h,f,g\), the finite
identities removing modulus primes likewise preserve a uniform
polynomial bound in \(N(a)\); the same interchange is thus valid
for each fixed sieve term in a sufficiently far right half plane.

It follows from \eqref{cub:lambda-pairing} that the full double
residue is
\begin{equation}\label{cub:complete-residue}
 \begin{split}
 \mathcal R_3(u)
 &=\zeta_\lambda(2)^2
   \sum_{a\equiv1\ (3)}\frac1{N(a)^u}
      \sum_{E,E'\in\mathcal E}
       H(E'E^{-1};u)v_E(a)\overline{v_{E'}(a)}\\
 &=\zeta_\lambda(2)^2d_\lambda(u)T_3(u).
 \end{split}
\end{equation}
This identity is first obtained in the region of uniform absolute
convergence just described. Formula \eqref{cub:theta-square-sum}
and meromorphic continuation then extend it to a neighborhood of
\(u_0\). Thus the diagonal square sum appears only after evaluation
of the full finite transformation. Using \(T_3(u_0)=0\), we obtain
\begin{equation}\label{cub:complete-residue-derivative}
 \mathcal R_3(u_0)=0,\qquad
 \mathcal R_3'(u_0)
 =\zeta_\lambda(2)^2d_\lambda(u_0)T_3'(u_0)>0.
\end{equation}
To pass from this unsieved residue to the squarefree family, we
must show that the sieve shifts preserve the finite pairing and
that the sieve sum may be interchanged with the Laurent coefficient.
The following translation formula will be used for the first step.

\begin{lemma}[Translation of the class residues]
\label{cub:class-translation-lemma}
Let \(\alpha\) be primary and squarefree, with \((\alpha,a)=1\).
Then
\begin{equation}\label{cub:class-translation}
 v_E(\alpha a)
 =N(\alpha)^{-\frac23}\overline{g_3(1,\alpha)}
      \chi_\alpha(a)\,v_{E[\alpha]}(a)
 \qquad(E\in\mathcal E).
\end{equation}
\end{lemma}

\begin{proof}
For the Fourier character \(\xi=\xi_{v,t}\), put
\(t_\xi=\omega^{2v}\lambda^{2t}\). Thus
\(\xi(c)=\overline{\left(\frac{t_\xi}{c}\right)_3}\).
Patterson's formula, \eqref{cub:raw-gauss-product}, and
\eqref{cub:raw-lambda-shifts} give
\begin{equation}\label{cub:shifted-fourier-residue}
 \rho(t_\xi\alpha a)
 =N(\alpha)^{-\frac23}\overline{g_3(1,\alpha)}
      \chi_\alpha(a)\,
      \overline{\xi(\alpha)}\,\rho(t_\xi a).
\end{equation}
Adjoining \(\alpha\) to the squarefree part at primes outside \(S\)
produces the first three factors. On the nonzero support with a
\(\lambda^2\) shift, the supplementary law contributes
\[
 \left(\frac{t_\xi}{\alpha}\right)_3
 =\overline{\xi(\alpha)}.
\]
On the remaining Fourier support both sides vanish.
By \eqref{cub:fourier-as-shift}, equation
\eqref{cub:shifted-fourier-residue} relates
\(\mathcal R_\xi(\alpha a)\) to \(\mathcal R_\xi(a)\).
Fourier inversion and
\[
 \overline{\xi(E)}\,\overline{\xi(\alpha)}
 =\overline{\xi(E[\alpha])}
\]
now prove \eqref{cub:class-translation}.
\end{proof}

\subsubsection{The finite decomposition of the squarefree sieve}
\label{cub:sieve-decomposition}

We continue with the natural family of Theorem~\ref{thm:second-moment}
and write
\[
 u_0=\frac16,\qquad w_0=\frac56.
\]
Every ideal coprime to \(3\) is represented by its unique primary generator.
Cubic reciprocity gives \(\chi_a(b)=\chi_b(a)\) for coprime primary
parameters. Thus the factor \(\psi_G(m_{[g]})\) in
\cite[Proposition~6.3.2]{DIPP} equals \(1\), and character orthogonality
simplifies the finite class projection in that proposition.

For a squarefree \(h\in I(S)\), put
\(S_h=\{p:p\mid h\}\) and \(T=S\cup S_h\).
We recall the precise series to which the sieve is applied.
For an ideal \(a\), let \(a_0\) be its part free of cube factors and
put \(a_1=\operatorname{rad}(a_0)\).
Let \(Q_a(s_1,s_2;\psi_1,\psi_2)\) be the correction polynomial of
\cite[Section~5.6.1]{DIPP}. In particular, \(Q_a=1\) when \(a\)
is squarefree. With their normalization, the perfect series is
\begin{equation}
\label{cub:sieve-perfect-definition}
 \begin{split}
 &\widetilde Z^T(s_1,s_2,w;\psi_1,\psi_2,\psi_3)\\
 &\quad=
 \sum_{a\in I(T)}
 \frac{\psi_3(a)Q_a(s_1,s_2;\psi_1,\psi_2)}{N(a)^w}
 L^T(s_1,\psi_1\chi_{a_0})
 L^T(s_2,\overline{\psi_2\chi_{a_0}}).
 \end{split}
\end{equation}
The conjugation here acts on the character, not on \(s_2\).
The finite sieve term of \cite[Sections~6.1.1--6.1.2]{DIPP} is
\begin{equation}
\label{cub:sieve-finite-series-definition}
 \begin{split}
 Z_1(s_1,s_2,w;h)
 ={}&
 \sum_{\substack{a\in I(S)\\h\mid\frac{a}{a_1}}}
 \frac{Q_a(s_1,s_2;1,1)}{N(a)^w}\\
 &\hspace{10mm}\times
 L^S(s_1,\chi_{a_0})L^S(s_2,\overline{\chi_{a_0}}).
 \end{split}
\end{equation}
For squarefree \(h\), the divisibility condition is equivalent to
\(\ord_p a>1\) for every \(p\mid h\).
M\"obius inversion therefore removes precisely the nonsquarefree
ideals. Since \(Q_a=1\) on the surviving ideals, it gives
\begin{equation}
\label{cub:sieve-mobius-identity}
 Z_1^S(s,s,w)=\sum_h^*\mu_K(h)Z_1(s,s,w;h)
\end{equation}
in the initial common domain of absolute convergence,
where the star restricts \(h\) to squarefree ideals in \(I(S)\).
We shall justify its use at \(w=w_0\) after removing the double pole.
Write
\[
 \widetilde Z_{f,g}^{\,T}(s,w)
 :=\widetilde Z^T(s,s,w;\chi_f,\chi_f,\chi_g),
 \qquad f,g\mid h^2,\quad(f,g)=1.
\]
We retain the convention of DIPP that the second input character
is conjugated in the second \(L\)-function. Their finite sieve
decomposition becomes
\begin{equation}
\label{cub:sieve-finite-decomposition}
 Z_1(s,s,w;h)
 =
 \sum_{\substack{f,g\mid h^2\\(f,g)=1}}
 f_f^{(1)}(s,s,w)\,g_{g,h(f)}(s,s,w)\,
 \chi_f(g)\,\widetilde Z_{f,g}^{\,T}(s,w),
\end{equation}
where \(h(f)\) is the part of \(h\) coprime to \(f\).
At a prime \(p\mid h\), there are exactly five possibilities:
\[
 (f_p,g_p)=(1,0),(2,0),(0,0),(0,1),(0,2),
 \qquad f_p=\ord_p f,\quad g_p=\ord_p g.
\]
We first identify the finite pairing in each term and cancel the
phases between different primes. We then compute the local Euler
weights and the five sieve coefficients.

\subsubsection{The finite pairing at \texorpdfstring{\(\lambda\)}{lambda} in the sieve}
\label{cub:sieve-lambda-stability}

When removing the finite Euler factors,
\cite[Proposition~5.4.4]{DIPP} twists both input characters by
\(\overline{\chi_t}\chi_{t'}\).
Here \(t,t'\) involve only the primes in \(T\setminus S\);
the factor at \(\lambda\) has already been included in \(d_\lambda(u)\).
The resulting twist is therefore
\[
 m=f\,\frac{t'}{t}.
\]
We allow \(m\) to be a fractional ideal and read its prime exponents
modulo \(3\). At a prime dividing \(g\), the outer character is
ramified, so the Euler sum has \((tt',g)=1\).
Together with \((f,g)=1\), this gives \((m,g)=1\).
Write
\[
 g=G_1G_2^2,\qquad
 G_1,G_2\ \text{squarefree and coprime},\qquad
 g_{\mathrm{rad}}=G_1G_2.
\]
After retaining the corrected root number, the two raw residue
indices are
\[
 a(mg_{\mathrm{rad}}g)^2,\qquad
 a(mg_{\mathrm{rad}}g^{-1})^2
 \quad\text{modulo cubes}.
\]
The outer index \(a\) belongs to \(I(T)\), so \((a,h)=1\).
For \(p\nmid g\), if \(m_p\equiv1\pmod3\), both indices contain a
factor \(p^2\) modulo cubes. Patterson's formula makes both residues
zero, including all their finite shifts at \(\lambda\).
For the remaining choices, set
\begin{equation}
\label{cub:sieve-surviving-indices}
 \delta=\prod_{\substack{p\mid h,\ p\nmid g\\m_p\equiv2\;(\mathrm{mod}\,3)}}p,
 \qquad B_L=\delta G_1,\qquad B_R=\delta G_2.
\end{equation}
The periodicity \eqref{cub:rho-periodicity} then reduces the two
indices to \(B_La\) and \(B_Ra\), respectively; the required
residue relations are those of \cite[Section~6]{DDFDS}.

\begin{lemma}[Translation of the finite pairing]
\label{cub:sieve-stable-pairing}
For every term of \eqref{cub:sieve-finite-decomposition} and every
nonzero choice of the Euler parameters, the finite pairing at
\(\lambda\) has the same scalar factor \(d_\lambda(u)\) as in the
unsieved case. The pairing includes the root number in
\eqref{cub:corrected-root-number}, all finite \(\gamma\)-factors,
and the Euler factors of the actual primitive characters.
\end{lemma}

\begin{proof}
The class translation formula \eqref{cub:class-translation} gives
\[
 v_E(B_La)
 =
 N(B_L)^{-\frac23}\overline{g_3(1,B_L)}
 \chi_{B_L}(a)\,v_{EB_L}(a),
\]
and the second residue satisfies the analogous identity.
After extracting the factors independent of \(E,E'\), change variables by
\[
 E_0=EB_L,\qquad E'_0=E'B_R.
\]
Since
\(\frac{B_R}{B_L}=\frac{G_2}{G_1}\equiv g^{-1}\pmod{K^{\times3}}\),
we have
\begin{equation}
\label{cub:sieve-class-difference}
 E'E^{-1}g^{-1}=E'_0E_0^{-1}.
\end{equation}
The finite conductor and root number in the corrected functional
equation depend on the difference class \([E'E^{-1}g^{-1}]\);
the actual Euler character at \(\lambda\) corresponds to its inverse.
Thus \eqref{cub:sieve-class-difference} transforms both the finite
root number kernel and the Euler ratio at \(\lambda\) into the
unsieved kernel \eqref{cub:lambda-kernel}.

We use \eqref{cub:gamma-cancellation} for the complete product of
\(\gamma\)-factors. After this cancellation,
\(C_{\mathrm{printed}}^{-1}\) contributes the \(g\)-phase
\(f_\lambda([g])^{-1}\), whereas
\(G(\overline{\chi_g})^{-1}\) contributes \(f_\lambda([g])\)
at \(\lambda\). Their product is \(1\).
The change of finite classes is a permutation, so it introduces
no additional normalization. By \eqref{cub:lambda-pairing}, the
remaining class sum is
\(d_\lambda(u)|\rho(a)|^2\), with
\[
 d_\lambda(u_0)=\frac{4}{9(1-3^{-\frac12})}>0.
\]
There is no remaining finite character of \(f,g,t,t'\), and no
additional factor equal to the number of classes.
\end{proof}

\subsubsection{Cancellation of phases between the sieve primes}
\label{cub:sieve-phase-cancellation}

We next combine the Gauss factors extracted in the preceding lemma
with \(\chi_f(g)\) and the phases from the Euler sum.
This step is needed to prove multiplicativity: real weights at
individual primes alone do not give an Euler product for the sieve.

\begin{lemma}[Cancellation between distinct primes]
\label{cub:sieve-cross-phase}
For every nonzero choice in \eqref{cub:sieve-surviving-indices},
the phases from the raw residue pairing, the corrected global root
number, the sieve coefficient \(\chi_f(g)\), and the Euler sum
have product \(1\). After extracting the finite pairing of
Lemma~\ref{cub:sieve-stable-pairing}, the remaining arithmetic
weight is
\begin{equation}
\label{cub:sieve-state-weight}
 N(\delta)^{-\frac13}N(g_{\mathrm{rad}})^{u-\frac23}.
\end{equation}
\end{lemma}

\begin{proof}
All products of Gauss sums below have coprime moduli.
Cubic reciprocity and twisted multiplicativity give
\[
 g_3(1,xy)=\overline{\chi_x(y)}\,g_3(1,x)g_3(1,y),
 \qquad
 |g_3(1,\delta)|^2=N(\delta).
\]
Consequently,
\begin{equation}
\label{cub:sieve-gauss-cross}
 \overline{g_3(1,\delta G_1)}g_3(1,\delta G_2)
 =
 \chi_\delta(g)N(\delta)\,
 \overline{g_3(1,G_1)}g_3(1,G_2).
\end{equation}
The ratio of the two residue indices satisfies
\(\frac{B_L}{B_R}\equiv g\pmod{K^{\times3}}\).
The resulting character of \(a\) is therefore canceled by the
outer factor \(\overline{\chi_g(a)}\).
The coefficient extracted from the raw pairing is
\begin{equation}
\label{cub:sieve-raw-mixed-coefficient}
 \chi_\delta(g)N(\delta)^{-\frac13}
 N(g_{\mathrm{rad}})^{-\frac23}
 \overline{g_3(1,G_1)}g_3(1,G_2).
\end{equation}

By \eqref{cub:general-g-root}, the cross factor between \(G_1\)
and \(G_2\) in the Chinese remainder decomposition of the primitive
root number \(G(\overline{\chi_g})\) is
\[
 \overline{\chi_{G_1}(G_2)}\chi_{G_2}(G_1)=1.
\]
Its part at the primes outside \(S\) is therefore
\[
 \frac{\overline{g_3(1,G_1)}g_3(1,G_2)}
 {N(g_{\mathrm{rad}})^{\frac12}}.
\]
After the finite phases have been canceled as in
Lemma~\ref{cub:sieve-stable-pairing}, division by this root number
changes the power \(N(g_{\mathrm{rad}})^{-\frac23}\) in
\eqref{cub:sieve-raw-mixed-coefficient} to
\(N(g_{\mathrm{rad}})^{-\frac16}\).
The factor \(N(g_{\mathrm{rad}})^{u-\frac12}\) outside the functional
equation then gives \(N(g_{\mathrm{rad}})^{u-\frac23}\).

The character coefficient from the Euler sum is
\(\chi_g(\frac{t'}{t})\), whereas the sieve coefficient in
\eqref{cub:sieve-finite-decomposition} is \(\chi_f(g)\), without
conjugation. Reciprocity identifies their product with \(\chi_m(g)\).
Combining this with the phase in
\eqref{cub:sieve-raw-mixed-coefficient}, we obtain
\[
 \chi_m(g)\chi_\delta(g)=\chi_{m\delta}(g)=1.
\]
Indeed, for every surviving choice the exponent
\(\ord_p(m\delta)\) is \(0\) modulo \(3\): its possibilities are
\(0\) and \(2+1\). This argument also covers negative exponents
of \(m\). Only the weight \eqref{cub:sieve-state-weight} remains.
\end{proof}

\subsubsection{Local weights from removing Euler factors}
\label{cub:sieve-cleanup}

Fix \(p\mid h\) and write \(q=Np\).
If \(p\nmid g\), the local sum in
\cite[Proposition~5.4.4]{DIPP}, after canceling the preceding
phases, is
\[
 \frac{1}{1-q^{-3u}}
 \sum_{j=0}^{2}\sum_{k=0}^{1}
 (-1)^kq^{-ju+k(u-1)}X^{k-j},
\]
Here \(X\) is the local operator that twists both input characters
by \(\chi_p\), so \(X^3=I\).
The exponent \(k-j\) records the conjugation in
\(\overline{\chi_t}\chi_{t'}\).

\begin{lemma}[The unramified and ramified weights]
\label{cub:sieve-local-weights}
For \(p\nmid g\), the local operator is
\begin{align}
\label{cub:sieve-cleanup-operator}
 K_q(u)
 &=(I-q^{u-1}X)(I-q^{-u}X^{-1})^{-1}\\
 \notag
 &=
 \frac{
 (1-q^{-1})I+(q^{-2u}-q^{u-1})X
 +q^{-u}(1-q^{-1})X^2
 }{1-q^{-3u}}.
\end{align}
Let \(v=(1,0,q^{-\frac13})\), with
\((Xv)_e=v_{e+1}\) and subscripts read modulo \(3\).
If \(k_e(q,u)=(K_q(u)v)_e\), then
\begin{equation}
\label{cub:sieve-k-values}
 \begin{aligned}
 k_0(q,u_0)&=(1+q^{-\frac12})^2,\\
 k_1(q,u_0)&=q^{-\frac16}+2q^{-\frac23},\\
 k_2(q,u_0)&=2q^{-\frac13}+q^{-\frac56}.
 \end{aligned}
\end{equation}
For \(p\mid g\), there is no Euler factor to remove at \(p\).
Both possibilities \(\ord_p g=1,2\) have weight
\begin{equation}
\label{cub:sieve-ramified-weight}
 q^{u-\frac23},
 \qquad\text{which equals }q^{-\frac12}\text{ at }u=u_0.
\end{equation}
\end{lemma}

\begin{proof}
Summing the finite expression above using \(X^3=I\) gives
\eqref{cub:sieve-cleanup-operator}.
When the exponent of the net twist at \(p\) is \(0,1,2\),
the exponents in both raw residue indices are \(0,2,1\),
respectively. Patterson's formula therefore gives the vector
of weights \(v=(1,0,q^{-\frac13})\).

Put \(b=q^{-\frac16}\). Then \(q^{-\frac13}=b^2\), and the three
numerators obtained by applying \eqref{cub:sieve-cleanup-operator}
to \(v\) are
\[
 (1-b^6)(1+b^3),\qquad
 b+b^4-2b^7,\qquad
 2b^2-b^5-b^8.
\]
Division by \(1-b^3\) gives
\((1+b^3)^2\), \(b+2b^4\), and \(2b^2+b^5\).
This proves \eqref{cub:sieve-k-values}.
The ramified weight follows directly from
\eqref{cub:sieve-state-weight}.
\end{proof}

\subsubsection{The five local sieve coefficients}
\label{cub:sieve-five-functions}

We calculate the sieve coefficients from the generating functions
for \(P_k,Q_k\) in \cite[Sections~5.5--5.6 and~6.2]{DIPP}.
For fixed \(q\), put
\[
 \begin{gathered}
 x=q^{-\frac12},\qquad y=q^{-w},\qquad Y=y^3,\\
 \Delta=(1-Y)^2\left(1-\frac{Y}{x}\right)^2,\qquad
 \mathcal E=(1-Y)(1-x^2Y)-\Delta.
 \end{gathered}
\]
Let \(F_1,F_2\) be the values of \(f_1^{(1)},f_1^{(2)}\)
at \((x,x,y;q)\), and let \(G_e\) be the value of \(g^{(e)}\)
there. These functions depend on \(q,w\).

\begin{lemma}[The local rational functions]
\label{cub:sieve-five-rational}
The following identities hold as convergent power series for
sufficiently small \(y\), and hence as identities of rational functions:
\begin{equation}
\label{cub:sieve-five-formulas}
 \begin{aligned}
 F_1&=y\left(\frac{1+2Y}{\Delta}-1\right),
 &
 F_2&=\frac{y^2(2+Y)}{\Delta},\\
 G_0&=\frac{3Y}{\Delta}
 +\frac{(1-x+x^2)\mathcal E}
 {\Delta(1-x)(1-x^3)},
 &
 G_1=G_2&=
 \frac{x\mathcal E}{\Delta(1-x)(1-x^3)}.
 \end{aligned}
\end{equation}
\end{lemma}

\begin{proof}
For \(\xi\in\mu_3\), set
\[
 B_\xi=(1-\xi x)(1-\xi^{-1}x).
\]
Since \(q(\xi x)(\xi^{-1}x)=1\), the polynomials of
\cite[Proposition~5.5.3]{DIPP} specialize to
\[
 P_0=1,\qquad
 P_{3n+e}(\xi x,\xi^{-1}x;q)=3n+e
 \quad(e=1,2,\ n\geq0),
\]
and
\[
 P_{3n}(\xi x,\xi^{-1}x;q)=3nB_\xi+1-x^2
 \quad(n\geq1).
\]
Their normalized generating function is
\[
 \sum_{k\geq0}Q_k(\xi x,\xi^{-1}x;q)y^k
 =
 \left(1-\frac{Y}{x}\right)^{-2}
 \sum_{k\geq0}P_k(\xi x,\xi^{-1}x;q)y^k.
\]
The normalizing factor contains only powers of \(y^3\), so it
preserves each residue class of exponents modulo \(3\).
Using
\[
 \sum_{n\geq0}(3n+1)Y^n=\frac{1+2Y}{(1-Y)^2},
 \qquad
 \sum_{n\geq0}(3n+2)Y^n=\frac{2+Y}{(1-Y)^2},
\]
we obtain
\[
 \sum_{n\geq0}Q_{3n+1}y^{3n+1}
 =\frac{y(1+2Y)}{\Delta},\qquad
 \sum_{n\geq0}Q_{3n+2}y^{3n+2}
 =\frac{y^2(2+Y)}{\Delta}.
\]
The squarefree sieve retains local exponents greater than \(1\).
We therefore subtract \(Q_1y=y\) from the first expression and
leave the second unchanged. This gives \(F_1,F_2\).

For exponents divisible by \(3\), the same calculation gives
\[
 \sum_{n\geq0}Q_{3n}y^{3n}
 =
 \frac{1+(3B_\xi-1-x^2)Y+x^2Y^2}{\Delta}.
\]
Subtracting the constant term and dividing by \(B_\xi\), we find
\begin{equation}
\label{cub:sieve-g-before-fourier}
 g(\xi x,\xi^{-1}x,y;q)
 =
 \frac{3Y}{\Delta}+\frac{\mathcal E}{B_\xi\Delta}.
\end{equation}
Fourier projection according to the difference between the
exponents of \(x_1,x_2\), modulo \(3\), gives
\[
 G_e=\frac13\sum_{j=0}^{2}
 \omega^{-ej}g(\omega^jx,\omega^{-j}x,y;q).
\]
Here \(B_1=(1-x)^2\) and
\(B_\omega=B_{\omega^2}=1+x+x^2\), so
\[
 \frac13\sum_{j=0}^{2}\frac1{B_{\omega^j}}
 =\frac{1-x+x^2}{(1-x)(1-x^3)}.
\]
For \(e=1,2\), we similarly have
\[
 \frac13\sum_{j=0}^{2}\frac{\omega^{-ej}}{B_{\omega^j}}
 =\frac{x}{(1-x)(1-x^3)}.
\]
Substituting these identities into
\eqref{cub:sieve-g-before-fourier} proves the remaining formulas.
\end{proof}

\subsubsection{Normalization at the intersecting poles}
\label{cub:sieve-laurent-normalization}

Write
\[
 \kappa_T=\Res_{s=1}\zeta_K^T(s),
 \qquad
 G_\infty(u)=
 \frac{\zeta_{K,\infty}(u)}{\zeta_{K,\infty}(1-u)}.
\]
For fixed \(h,f,g\), let \(\mathcal R_{f,g}^{\,T}(u)\) denote
the complete finite pairing obtained by taking the residues in
both Kubota variables in the corrected functional equation,
including the factor \(N(g_{\mathrm{rad}})^{u-\frac12}\).
Its normalization is fixed by
\begin{equation}
\label{cub:sieve-double-residue-identity}
 \begin{split}
 &\lim_{s\to u+\frac13}\left(s-u-\frac13\right)^2
       \widetilde Z_{f,g}^{\,T}(s,1-u)\\
 &\hspace{15mm}
 =\zeta_K^T\left(2u+\frac23\right)G_\infty(u)
   \mathcal R_{f,g}^{\,T}(u).
 \end{split}
\end{equation}
This follows initially from the functional equation in
\cite[Section~5.7]{DIPP}, with the root number
\eqref{cub:corrected-root-number}, in a right half-plane where
the outer residue series is absolutely convergent.
The Cauchy integral argument used for \eqref{cub:complete-residue}
justifies taking the residues term by term.
The identity then continues meromorphically.

Define
\[
 T_3^h(u)=
 \sum_{\substack{a\in I(S)\\(a,h)=1}}
 \frac{|\rho(a)|^2}{N(a)^u},
 \qquad
 c_h=\prod_{p\mid h}\bigl(1-(Np)^{-1}\bigr),
\]
and
\[
 \mathcal W_{f,g}(u)=
 \prod_{\substack{p\mid h\\p\nmid g}}
       k_{\ord_p f}(Np,u)
 \prod_{p\mid g}(Np)^{u-\frac23}.
\]
Lemmas~\ref{cub:sieve-stable-pairing},
\ref{cub:sieve-cross-phase}, and~\ref{cub:sieve-local-weights},
together with the removal formula
\eqref{cub:normalized-residue} for each normalized residue, give
\begin{equation}
\label{cub:sieve-complete-residue-pairing}
 \chi_f(g)\mathcal R_{f,g}^{\,T}(u)
 =
 \zeta_K^S(2)^2c_h^2d_\lambda(u)
 T_3^h(u)\mathcal W_{f,g}(u).
\end{equation}
The factor \(\chi_f(g)\) remains on the left.
It is precisely the sieve coefficient needed to cancel the
phases between distinct primes.

\begin{lemma}[The coefficient of the double pole]
\label{cub:sieve-normalization-lemma}
Let
\[
 A_{\mathrm{perf}}
 =
 \lim_{w\to w_0}(w-w_0)^2
       \widetilde Z^S\left(\frac12,\frac12,w;1,1,1\right).
\]
Then
\begin{equation}
\label{cub:sieve-Aperf}
 \begin{split}
 A_{\mathrm{perf}}
 &=\frac{\kappa_S}{2}G_\infty(u_0)
     \zeta_K^S(2)^2d_\lambda(u_0)T_3'(u_0)\\
 &=(2\pi)^{\frac23}
   \frac{\Gamma\left(\frac16\right)}{\Gamma\left(\frac56\right)}
   \zeta_K^S(2)^2d_\lambda\left(\frac16\right)
   \rho_0^2\zeta_K^S\left(\frac12\right)^2>0.
 \end{split}
\end{equation}
For \(T=S\cup S_h\), the common local normalization at each
sieve prime \(p\mid h\), relative to \(A_{\mathrm{perf}}\), is
\begin{equation}
\label{cub:sieve-Dp}
 D_p=(1-q^{-1})^3
       \frac{1-q^{-\frac12}}{1+q^{-\frac12}},
 \qquad q=Np.
\end{equation}
\end{lemma}

\begin{proof}
By \cite[Section~7.0.2]{DIPP},
\[
 H_{f,g}^{\,T}(s,w)
 :=\left(s+w-\frac43\right)^2\widetilde Z_{f,g}^{\,T}(s,w)
\]
is holomorphic near \((\frac12,\frac56)\).
On the polar line, put \(s=u+\frac13\) and \(w=1-u\).
Equation~\eqref{cub:sieve-double-residue-identity} gives
\[
 H_{f,g}^{\,T}\left(u+\frac13,1-u\right)
 =\zeta_K^T\left(2u+\frac23\right)G_\infty(u)
       \mathcal R_{f,g}^{\,T}(u).
\]
Equations~\eqref{cub:sieve-complete-residue-pairing} and
\eqref{cub:theta-derivative} show that
\(\mathcal R_{f,g}^{\,T}\) is holomorphic at \(u_0\) and vanishes there.
Since
\[
 \zeta_K^T\left(2u+\frac23\right)
 =\frac{\kappa_T}{2(u-u_0)}+O(1),
\]
we obtain
\begin{equation}
\label{cub:sieve-collision-value}
 \lim_{w\to w_0}(w-w_0)^2
 \widetilde Z_{f,g}^{\,T}\left(\frac12,w\right)
 =
 \frac{\kappa_T}{2}G_\infty(u_0)
       (\mathcal R_{f,g}^{\,T})'(u_0).
\end{equation}
On the slice \(s=\frac12\), the factor removing the pole is
exactly \((w-w_0)^2\), so there is no further change of scale or sign.
Taking \(h=f=g=1\) and substituting
\eqref{cub:theta-derivative} proves \eqref{cub:sieve-Aperf}.

For general \(h\), the local factors of
\eqref{cub:theta-square-sum} give
\begin{equation}
\label{cub:sieve-restricted-theta}
 T_3^h(u)=T_3(u)
 \prod_{p\mid h}\frac{1-q^{-3u}}{1+q^{-u-\frac13}}.
\end{equation}
As \(T_3(u_0)=0\), this implies
\[
 \frac{(T_3^h)'(u_0)}{T_3'(u_0)}
 =
 \prod_{p\mid h}
 \frac{1-q^{-\frac12}}{1+q^{-\frac12}}.
\]
Moreover,
\[
 \frac{\kappa_T}{\kappa_S}
 =\prod_{p\mid h}(1-q^{-1}),
 \qquad
 c_h^2=\prod_{p\mid h}(1-q^{-1})^2.
\]
Multiplying these three contributions gives \eqref{cub:sieve-Dp}.
Two of the three factors \(1-q^{-1}\) come from the normalized
Kubota residues, and the third comes from the outer zeta residue.

When \eqref{cub:sieve-complete-residue-pairing} is differentiated
at \(u_0\), the derivatives of \(d_\lambda\), \(\mathcal W_{f,g}\),
and the local correction in \eqref{cub:sieve-restricted-theta}
are all multiplied by \(T_3(u_0)=0\).
Their contributions therefore vanish.
This uses the common zero of the complete finite pairing after
that pairing has been evaluated.
\end{proof}

\begin{proposition}[The contribution of a sieve prime]
\label{cub:sieve-local-polynomial}
Put
\[
 A_h=\lim_{w\to w_0}(w-w_0)^2
 Z_1\left(\frac12,\frac12,w;h\right).
\]
For every squarefree \(h\),
\begin{equation}
\label{cub:sieve-Ah-factorization}
 A_h=A_{\mathrm{perf}}\prod_{p\mid h}R_p,
\end{equation}
where, with \(q=Np\),
\begin{equation}
\label{cub:sieve-Rp}
 R_p=
 6q^{-2}-2q^{-\frac52}-8q^{-3}
 +4q^{-\frac72}+3q^{-4}-2q^{-\frac92}.
\end{equation}
Equivalently,
\begin{equation}
\label{cub:sieve-positive-local-factor}
 1-R_p
 =(1-q^{-1})^2
       (1+2q^{-1}-3q^{-2}+2q^{-\frac52})>0.
\end{equation}
\end{proposition}

\begin{proof}
Equations~\eqref{cub:sieve-finite-decomposition} and
\eqref{cub:sieve-complete-residue-pairing} show that, after taking
the coefficient of the double pole, the choices at distinct primes
can be multiplied independently.
For a fixed \(p\mid h\), the contribution is
\begin{equation}
\label{cub:sieve-five-branch-combination}
 R_p
 =
 D_p\left[
 F_1k_1+F_2k_2+G_0k_0+(G_1+G_2)q^{-\frac12}
 \right]_{w=w_0,\ u=u_0}.
\end{equation}
The five sieve functions are holomorphic at this point, so only
their values enter the coefficient of the double pole.
Their derivatives in \(w\) can contribute only to the coefficient
of the simple pole.

To simplify \eqref{cub:sieve-five-branch-combination}, put
\(z=q^{-\frac12}\). At the point in question,
\[
 Y=z^5,\qquad
 \Delta=(1-z^5)^2(1-z^4)^2,\qquad
 \mathcal E=(1-z^5)(1-z^7)-\Delta,\qquad
 D_p=(1-z)^4(1+z)^2.
\]
The first two terms are
\[
 F_1k_1
 =\frac{z^2(1+2z)(1+2z^5-\Delta)}{\Delta},
 \qquad
 F_2k_2
 =\frac{z^4(2+z)(2+z^5)}{\Delta}.
\]
For the remaining three terms, the coefficient of \(\mathcal E\)
simplifies by
\[
 \frac{(1-z+z^2)(1+z)^2+2z^2}
      {(1-z)(1-z^3)}
 =\frac{1+z^2}{(1-z)^2}.
\]
Writing
\[
 P_0=
 z^2(1+2z)(1+2z^5-\Delta)
 +z^4(2+z)(2+z^5)+3z^5(1+z)^2,
\]
we consequently obtain
\[
 R_p=
 \frac{(1-z)^4(1+z)^2P_0
       +(1-z)^2(1+z)^2\mathcal E(1+z^2)}
      {\Delta}.
\]
Substitution of \(\Delta,\mathcal E\) shows that its numerator is
\[
 \Delta\bigl(6z^4-2z^5-8z^6+4z^7+3z^8-2z^9\bigr).
\]
This proves \eqref{cub:sieve-Rp}.
The identity
\[
 1-6z^4+2z^5+8z^6-4z^7-3z^8+2z^9
 =(1-z^2)^2(1+2z^2-3z^4+2z^5)
\]
then gives \eqref{cub:sieve-positive-local-factor}.
The second factor is positive because
\[
 1+2q^{-1}-3q^{-2}+2q^{-\frac52}
 =(1-q^{-1})(1+3q^{-1})+2q^{-\frac52}>0.
\]
\end{proof}

\subsubsection{Normal convergence and the squarefree coefficient}
\label{cub:sieve-normal-convergence}

It remains to identify the sum of the coefficients \(A_h\) with
the Laurent coefficient of the actual squarefree series.
For this purpose we need normal convergence of the sieve after
removing the pole.

\begin{lemma}[Normal convergence of the regularized sieve]
\label{cub:sieve-normal-lemma}
There is \(\varrho>0\) such that
\begin{equation}
\label{cub:sieve-normal-sum}
 \sum_h^*
 \sup_{|w-w_0|\leq\varrho}
 \left|(w-w_0)^2
 Z_1\left(\frac12,\frac12,w;h\right)\right|<\infty.
\end{equation}
Consequently,
\begin{equation}
\label{cub:sieve-interchange}
 A_{\mathrm{sf}}
 :=
 \lim_{w\to w_0}(w-w_0)^2
 Z_1^S\left(\frac12,\frac12,w\right)
 =
 \sum_h^*\mu_K(h)A_h,
\end{equation}
and the series on the right is absolutely convergent.
\end{lemma}

\begin{proof}
Choose \(0<\varrho<\frac1{12}\) and put
\(\sigma_0=\frac56-\varrho>\frac34\).
Combining \cite[Propositions~6.3.2 and~7.2.2; Section~8.1.1]{DIPP},
we have, uniformly on this closed disc,
\begin{equation}
\label{cub:sieve-uniform-majorant}
 \left|(w-1)^2(w-w_0)^2
 Z_1\left(\frac12,\frac12,w;h\right)\right|
 \ll_{\varrho,\eps}
 N(h)^{2-4\sigma_0+\eps}M_h(\eps),
\end{equation}
where we may take
\[
 M_h(\eps)=
 \sum_{\psi\in\widehat R_{\mathfrak c}}
 \ \sum_{t\mid h^2}
 \ \sum_{\substack{a\in I(S),\ a\ \text{cube free}\\(a,t)=1}}
 \frac{\left|L\left(\frac12,\psi\chi_{at}\right)\right|^2}
 {N(a)^{1+\eps}}.
\]
Here we have used the uniform bound
\(O(q^{-2\Re w})\) for the five local sieve functions.
The number of finite choices and the other local constants are
absorbed into \(N(h)^\eps\).
The extra root number in \eqref{cub:corrected-root-number} has
absolute value \(1\), so it does not affect this estimate.
The actual Euler ratios at \(\lambda\) are fixed rational functions
on the nine classes. Their denominators do not vanish on the chosen
disc, and they are independent of \(h\); thus they change only the
implied constant in \eqref{cub:sieve-uniform-majorant}.

To prove summability, fix \(\sigma>\frac34\) and choose
\(0<\eps<4\sigma-3\).
In the nonnegative series
\[
 \sum_h^*N(h)^{2-4\sigma}M_h(\eps),
\]
write \(h=\operatorname{rad}(t)b\).
The sum over squarefree \(b\) converges because
\(2-4\sigma<-1\).
With \(\delta_0=4\sigma-3>0\), the remaining weight satisfies
\[
 \frac{1}
 {N(\operatorname{rad}t)^{1+\delta_0}N(a)^{1+\eps}}
 \leq
 \frac{1}{N(\operatorname{rad}(at))^{1+\eps}},
 \qquad(a,t)=1.
\]
After setting \(b_0=at\), the number of decompositions is at most
\[
 2^{\omega(b_0)}
 \ll_\eps N(\operatorname{rad}b_0)^{\frac{\eps}{2}}.
\]
The nonnegative series is therefore bounded by a constant times
\[
 \sum_{\psi\in\widehat R_{\mathfrak c}}
 \ \sum_{\substack{b_0\in I(S)\\b_0\ \text{cube free}}}
 \frac{\left|L\left(\frac12,\psi\chi_{b_0}\right)\right|^2}
 {N(\operatorname{rad}b_0)^{1+\frac{\eps}{2}}}.
\]
The second moment bound ordered by conductor in
\cite[Theorem~1.2 and Section~9, proof of Theorem~1.3]{deFaveri},
followed by partial summation, proves convergence of this last series.
Indeed, \(N(\operatorname{Cond}(\chi_{b_0}))\) differs from
\(N(\operatorname{rad}b_0)\) only by a factor from a fixed finite
set at \(\lambda\).

Now choose \(\eps\) in \eqref{cub:sieve-uniform-majorant} sufficiently
small that
\[
 \sigma=\sigma_0-\frac{\eps}{4}>\frac34,
 \qquad \eps<4\sigma-3.
\]
The right side of that estimate is then summable over \(h\).
Since \(w-1\) stays away from \(0\) on the disc,
\eqref{cub:sieve-normal-sum} follows.
Finally, the M\"obius identity \eqref{cub:sieve-mobius-identity}
and the continuation
in \cite[Section~9]{deFaveri} identify the normally convergent sum
with \((w-w_0)^2Z_1^S(\frac12,\frac12,w)\).
Evaluation at \(w=w_0\) gives \eqref{cub:sieve-interchange}.
\end{proof}

\begin{proposition}[Positivity of the squarefree coefficient]
\label{cub:sieve-final-coefficient}
For the natural primary cubic family of
Theorem~\ref{thm:second-moment}, with \(q=Np\),
\begin{equation}
\label{cub:sieve-final-euler}
 A_{\mathrm{sf}}
 =
 A_{\mathrm{perf}}
 \prod_{p\notin S}
 (1-q^{-1})^2
 \bigl(1+2q^{-1}-3q^{-2}+2q^{-\frac52}\bigr)
 >0,
\end{equation}
where \(A_{\mathrm{perf}}\) is given by \eqref{cub:sieve-Aperf}.
\end{proposition}

\begin{proof}
Substituting \eqref{cub:sieve-Ah-factorization} into the absolutely
convergent series \eqref{cub:sieve-interchange}, we obtain
\[
 A_{\mathrm{sf}}
 =A_{\mathrm{perf}}
 \sum_h^*\mu_K(h)\prod_{p\mid h}R_p
 =A_{\mathrm{perf}}\prod_{p\notin S}(1-R_p).
\]
By \eqref{cub:sieve-Rp},
\(R_p=6q^{-2}+O(q^{-\frac52})\), so this product converges absolutely.
Each factor is strictly positive by
\eqref{cub:sieve-positive-local-factor}, and
\(A_{\mathrm{perf}}>0\). This proves the assertion.
\end{proof}

\begin{remark}
\label{cub:sieve-product-scope}
The two factors at each prime in \eqref{cub:sieve-final-euler}
must be kept together. In particular,
\(\prod_p(1-q^{-1})^2\) is not a separate convergent positive product.
The calculation uses \(K=\Q(\omega)\), the fixed generator
\(\lambda=1+2\omega\), and all primary parameters.
Restricting the parameters to one class modulo \(9\), restoring
the Euler factor at \(\lambda\), or changing the base field requires
a new calculation of the finite projection and local pairing.
\end{remark}

We now complete the proof of Theorem~\ref{thm:second-moment}.
To express the coefficient in the normalization used there,
\eqref{cub:rho-constant} and \eqref{cub:lambda-pairing} give
\[
 \zeta_\lambda(2)^2\rho_0^2
 =\frac{(2\pi)^{\frac{10}{3}}}
 {3^7\Gamma\left(\frac23\right)^2},
 \qquad
 d_\lambda\left(\frac16\right)
 =\frac23\left(1+3^{-\frac12}\right).
\]
Since \(\zeta_\lambda=\zeta_K^S\), substituting these identities
into \eqref{cub:sieve-Aperf} and \eqref{cub:sieve-final-euler}
gives precisely \eqref{eq:second-coefficient}.
Let \(W\) be the nonnegative, nonzero smooth weight in that theorem,
and write
\[
 \widehat W(w)=\int_0^\infty W(x)x^{w-1}\,dx.
\]
If \(B_{\mathrm{sf}}\) denotes the coefficient of the simple pole, then
\[
 Z_1^S\left(\frac12,\frac12,w\right)
 =\frac{A_{\mathrm{sf}}}{(w-w_0)^2}
  +\frac{B_{\mathrm{sf}}}{w-w_0}+O(1).
\]
The defining series has real coefficients in its initial domain,
so its meromorphic continuation satisfies Schwarz symmetry.
Both Laurent coefficients at the real pole \(w_0\) are therefore real.
Mellin inversion and the contour shift of
\cite[Theorem~1.3]{deFaveri} give the secondary contribution
\begin{equation}
\label{cub:sieve-secondary-polynomial}
 \begin{split}
 X^{\frac56}\biggl\{
 &A_{\mathrm{sf}}\widehat W\left(\frac56\right)\log X\\
 &+A_{\mathrm{sf}}\widehat W'\left(\frac56\right)
  +B_{\mathrm{sf}}\widehat W\left(\frac56\right)
 \biggr\}.
 \end{split}
\end{equation}
Since \(\widehat W(\frac56)>0\), the coefficient of \(\log X\) is
strictly positive. Thus the secondary polynomial has degree exactly
\(1\), which proves the asserted nonvanishing.
No part of this second moment argument uses GRH.

\end{document}